\documentclass[colorinlistoftodos]{article}[11pt]
\usepackage{graphicx,hyperref,amsmath,amssymb,natbib,bm,url, amsthm}
\usepackage{microtype,todonotes}
\usepackage{subcaption}
\usepackage[australian]{babel}
\usepackage[a4paper,text={16.5cm,25.2cm},centering]{geometry}
\usepackage[compact,small]{titlesec}
\usepackage{enumitem}
\usepackage{tikz}
\usepackage{pgfplots}
\pgfplotsset{compat=1.18}
\usetikzlibrary{arrows,automata}
\usepackage{cleveref}
\usepackage{comment}
\usepackage{natbib}
\usepackage{wrapfig}

\usepackage{kpfonts}
\usepackage[T1]{fontenc}
\title{Random Unbiased Perturbations in Nonparametric Regression}
\author{Anna Lyubarskaja and Dominik Rothenh\"ausler}
\newtheorem{remark}{Remark}[section]
\newtheorem{lemma}{Lemma}[section]
\newtheorem{theorem}{Theorem}[section]

\newtheorem{prop}{Proposition}[section]
\newtheorem{assume}{Assumption}
\newtheorem{defn}{Definition}
\newtheorem{cor}{Corollary}[section]
\newcommand{\eff}{\mathrm{eff}}
\newcommand{\E}{\mathbb{E}}

\newcommand{\Be}{B_{\varepsilon}}
\newcommand{\e}{\varepsilon}

\newcommand{\Var}{\mathbb{V}\mathrm{ar}} 
\newcommand{\Cov}{\mathbb{C}\mathrm{ov}} 
 
\newcommand{\KL}{\mathrm{KL}}
\newcommand{\R}{\mathbb{R}}
\newcommand{\1}{\mathbf{1}}

\begin{document}
	\maketitle
\begin{abstract}
 We study nonparametric regression with covariates $X$ and outcome $Y$ under random unbiased perturbations (RUPs) of the conditional distribution $Y|X$, where the marginal distribution of covariates, $P^X$, remains fixed but the conditional law, $P^{Y|X}$, varies randomly across datasets. Unlike adversarial distribution shift frameworks that yield conservative worst-case guarantees, RUPs induce dataset-level variance inflation rather than systematic bias. We provide examples of RUPs and show that this distributional uncertainty reduces the effective sample size to $n_{\mathrm{eff}} = n/(1 + n \tau)$, where $\tau\in [0,1]$ quantifies the perturbation strength. For local polynomial estimators, we derive an extended bias-variance decomposition that includes a distributional variance term with the same bandwidth scaling as classical sampling variance. This leads to a modified bandwidth selection principle: when distributional uncertainty dominates sampling uncertainty ($\tau \gg 1/n$), optimal bandwidths scale as $\tau^{1/(2\beta+1)}$ rather than the usual $n^{-1/(2\beta+1)}$, where $\beta$ indicates the smoothness of the function class considered. We also establish matching minimax lower bounds showing that there exists an RUP for which this effective sample size $n_\eff$ is fundamental. Our results demonstrate that random dataset-level perturbations create a distinct mode of uncertainty that affects both practical tuning and fundamental statistical limits.

\end{abstract}

\section{Introduction}

Data are rarely sampled i.i.d. from the true distribution of interest. Instead, we often observe data from a slightly shifted version of our target distribution. For instance, the true cost of living for an average household is typically approximated by consumer price index (CPI), which is constructed by choosing a basket of representative goods and services. Data can then be collected to estimate this CPI, but the discrepancy between the true average cost of living and the CPI depends on the choice of goods and services included in the CPI basket, and is not reduced by collecting additional data. Similarly, clinical outcomes observed at a single hospital differ slightly from regional population-level outcomes as they depend on the hospital’s particularities, such as specific staff or protocols. As a final example, flight delays on a given day are influenced by that day’s weather patterns and will thus deviate from a long-run average distribution of delays. In expectation over many such days, hospitals or consumer baskets, one recovers a true target distribution $P_0$, but any single dataset is generated under just one realization $P_\xi$, a shifted version of the truth.

In the examples above, the distribution of covariates $X$ can be regarded as fixed (consumer spending patterns, patient populations, scheduled flights), while the conditional law of outcomes $Y\mid X$ is perturbed by measurement procedure or environment. We refer to such shifts as \emph{random unbiased perturbations} (RUPs): random, mean-zero deviations from a target law that preserve the covariate distribution. We define this notion formally in \Cref{def:rup-yx} (\Cref{sec:rups}).

Our RUP framework thus involves two distinct layers of randomness. First, $\Xi$ represents a distribution over possible perturbations, where a random realization $\xi \sim \Xi$ parameterizes the distributional shift. Second, conditional on the chosen $\xi$, the data $\{(X_i,Y_i)\}_{i=1}^n$ are sampled i.i.d.\ with $X_i \sim P^X$ as usual and $Y_i \mid X_i \sim P_\xi(\cdot \mid X_i)$. This perspective introduces a new source of uncertainty alongside the sampling error, as now the discrepancy between the observed empirical law, $\hat P_\xi$, and the target, $P_0$, naturally decomposes as
\begin{equation}
    \label{eq:err_decomp}
   \hat P_\xi - P_0 = (\hat P_\xi - P_\xi) + (P_\xi - P_0).
\end{equation}
The first term is the a familiar sampling error due to observing only $n$ samples from a given distribution, while the second term reflects a distributional shift from the perturbation $P_\xi \neq P_0$. 
\begin{figure}[t]
    \centering
        \includegraphics[width=\linewidth]{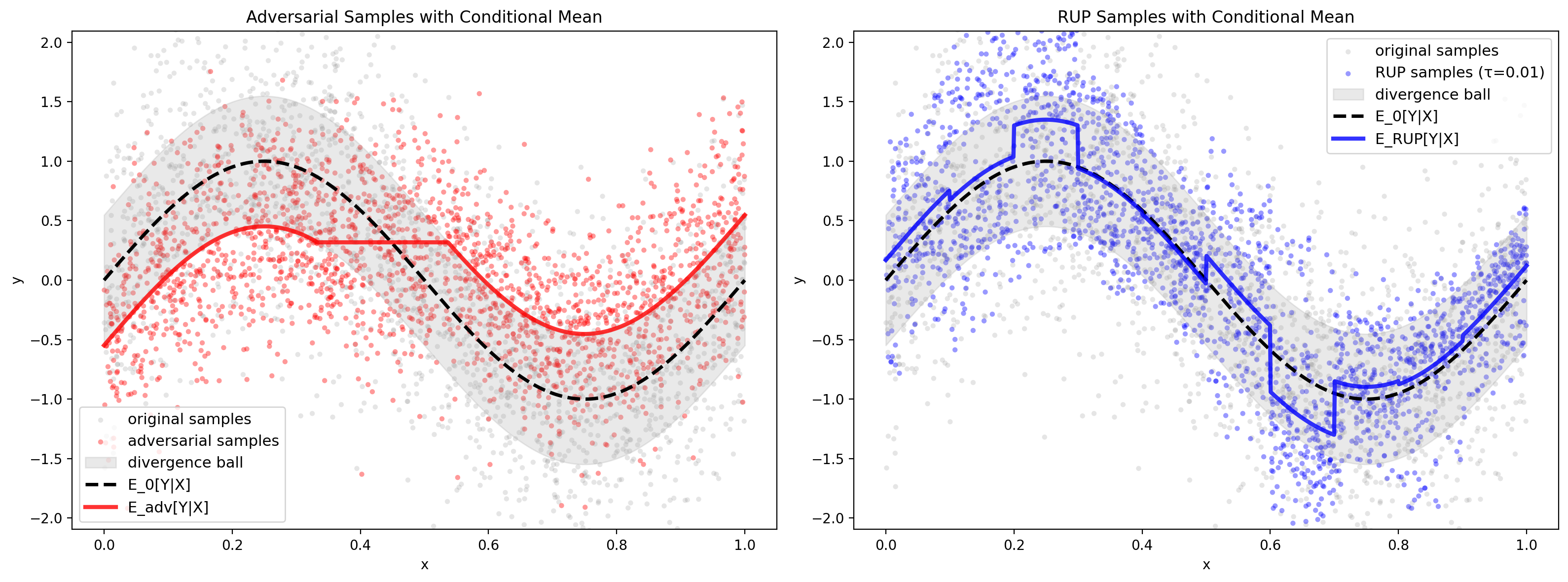}
        \caption{ Examples of points sampled adversarially (left, red) and from an RUP (right, blue). The two distribution shifts have the same KL divergence $\kappa$ from the baseline distribution $P_0$ (sampled in gray). However, the adversarial conditional mean (left, red) provides a much poorer approximation of the target $P_0$ conditional mean (black, dashed), whereas the RUP conditional mean (right, blue) remains much closer to the true relationship. The shaded gray divergence ball illustrates the range of conditional means that distributions at KL divergence $\kappa$ from $P_0$ can take.}
    \label{fig:why_no_dro}
\end{figure}

Much of the distribution shift literature relies on distributionally robust optimization, and treats the conditional $Y|X$ shift in the second term in \eqref{eq:err_decomp} adversarially. $P_\xi$ is assumed to lie within a divergence ball around $P_0$, and statistical methods are conservatively designed to be valid for all distributions within this ball\citep{duchi_learning_2021, cai_transfer_2021, reeve_adaptive_2021}. In contrast, we study perturbations that are random and mean-zero. \Cref{fig:why_no_dro} illustrates why the DRO approach is too conservative for the study of RUPs. We simulate two different $Y|X$ shifts of a baseline distribution $P_0$, each with the same KL divergence $\kappa$. On the left we visualize an adversarially shifted distribution, while the right side pictures a representative sample of an RUP. A DRO approach is optimized for adversarial shifts, which may distort the conditional mean systematically producing biased approximations of $\E_0[Y|X]$. RUPs, on the other hand, introduce only random, mean-zero fluctuations in statistics such as $\E[Y|X]$ that average out and stay near the target function. 

To build intuition for the behavior of RUPs, recall the flight delay example. Let $X$ represent the route information and $Y$ the delay. The marginal distribution of routes is fixed by the flight schedule, but a perturbation $\xi$ (i.e. a single day's weather) changes the conditional distribution of delays. These shifts are not uniform across all $X$: a snowstorm may delay flights in and out of Chicago, while flights elsewhere remain largely unaffected. Thus, RUPs may induce correlated changes in $Y\mid X$ across subsets of $X$, reflecting the geography of the covariates. We summarize the overall impact of an RUP by a perturbation strength parameter $\tau$, which we later decompose into a variance scale $\delta^2$ and an $X$-dependency parameter $\bar{\rho}$ (\Cref{def:rup-yx}). 

\paragraph{An extended bias–variance decomposition.} Under the RUP setting, the distributional shift term in \eqref{eq:err_decomp}, $P_\xi - P_0$, appears as an additional variance term in the error. Consider for instance the pointwise risk. Let $f^0(x) = \E_0[Y \mid X=x]$ denote the baseline regression function under $P_0$, and let $f^\xi(x) = \E_{\xi}[Y \mid X=x]$ denote the regression function under a perturbed law $P_\xi$, with estimator $\hat{f}^\xi_n$ constructed from $n$ samples drawn i.i.d. from $P_\xi$. We let the subscript $\Xi$ denote "with respect to the distribution over perturbations $\xi' \sim \Xi$", and the subscript $\xi$ denote "conditional on the realization $\xi$".

Applying the law of total variance to $\hat f_n^\xi(x_0)$ yields
\begin{equation}
\label{eq:bv_rup}
\E_\Xi \E_{\xi}\big[(\hat f_n^\xi(x_0)-f^0(x_0))^2\big]
= \underbrace{\big(\E_\Xi \E_{\xi}[\hat f_n^\xi(x_0)] - f^0(x_0)\big)^2}_{\text{bias squared}}
+ \underbrace{\E_\Xi\big[\Var_{\xi}(\hat f_n^\xi(x_0))\big]}_{\text{sampling variance}}
+ \underbrace{\Var_\Xi\big(\E_{ \xi}[\hat f_n^\xi(x_0)]\big)}_{\text{distributional variance}}.
\end{equation}

The first two terms coincide with the usual bias variance decomposition under i.i.d.\ sampling from $P_0$. 
The third term is new: it captures dataset-level variability introduced by the perturbation. 
\begin{figure}[t]
    \centering
\begin{tikzpicture}
\begin{axis}[
  width=0.7\linewidth, height=0.45\linewidth,
  xmin=0.2, xmax=11,
  ymin=-0.06, ymax=0.8,
  axis lines=left,
  axis x line=middle,
  tick style={black, thick},
  xtick={3.0571,4.6415},
  xticklabels={$c^\star_{\text{RUP}}$, $c^\star_{\text{classical}}$}, 
  x tick label style={anchor=north, yshift=-2pt}, 
  xlabel={model complexity},
  xlabel style={at={(axis description cs:0.5,-0.08)},anchor=north}, 
  ytick=\empty,
  ylabel={risk (MSE)},
  legend style={
    at={(0.5, 1.2)}, anchor=north,
    draw=black, fill=white, rounded corners,
    legend columns=3, font=\small
  },
  domain=0.2:11, samples=200, thick
]
\def\beta{1.0}
\def\ninv{0.02}
\def\minv{0.05}


\addplot[black, dashed] {(1/x)^(2*\beta)};                     \addlegendentry{bias$^2$}

\addplot[green!70!black, dashed] {\ninv*x};                    \addlegendentry{sampling variance}

\addplot[blue, thick] {(1/x)^(2*\beta) + \ninv*x};             \addlegendentry{MSE (no shift)}

\addlegendimage{empty legend}
\addlegendentry{}

\addplot[orange, dashed] {\minv*x};                            \addlegendentry{distributional variance}

\addplot[red,  thick] {(1/x)^(2*\beta) + \ninv*x + \minv*x};   \addlegendentry{MSE (RUP)}

\addplot[blue, mark=*, only marks]
  coordinates {(4.6415,{(1/4.6415)^(2*\beta) + \ninv*4.6415})};
\addplot[red,  mark=*, only marks]
  coordinates {(3.0571,{(1/3.0571)^(2*\beta) + \ninv*3.0571 + \minv*3.0571})};

\end{axis}
\end{tikzpicture}

    \caption{Schematic bias–variance tradeoff under random perturbations. 
Bias decreases with model complexity, while sampling variance and distributional variance both increase. The additional distributional variance raises the total RUP-adjusted curve and shifts the optimal complexity leftward (toward more regularized models).}
    \label{fig:bv_shift}
\end{figure}
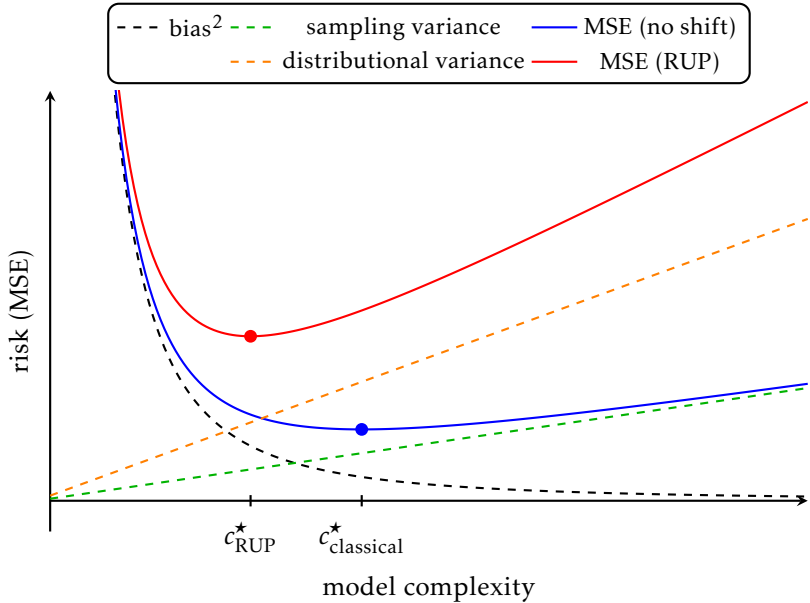
\paragraph{Our contributions.}
We study nonparametric regression under RUPs of the conditional law $Y\mid X$. 
A central consequence of RUPs is that they reduce the effective sample size available for estimation 
\begin{equation}
    \label{eq:n_eff}
    n_\eff = \frac{n}{1 + n\tau},
\end{equation}
where $\tau$ is a measure of strength of the RUP (as outlined in \Cref{def:rup-yx}).

\begin{enumerate}[leftmargin=*]

\item \textbf{Variance inflation and optimal tuning.}  
For local polynomial estimators (LPEs), perturbations introduce an additional variance term with the same complexity dependence as classical sampling variance. Ignoring this term leads to overfitting and underestimation of risk. Accounting for it yields a modified tuning principle: bandwidths should be chosen by replacing $n$ with $n_\eff$. An illustration of this adjusted bandwidth tuning is presented in \Cref{fig:bv_shift}. This interpretation highlights that under RUPs, distributional uncertainty acts like variance rather than bias.
\item \textbf{Fundamental limits.}  
Using minimax lower bounds, we show that $n_\eff$ also captures the fundamental information limit. We provide an example of an RUP set up where the amount of information scales with $n_\eff$ rather than with $n$, showing that distributional uncertainty fundamentally caps the rate of convergence of any estimator. In particular, we show that in the regime where the sampling and distributional variance are of the same order, i.e. the perturbation strength scales as $\tau_n \propto \frac{1}{n}$, $n_\eff$ captures the correct rescaling between the two for the optimal convergence rate.
\end{enumerate}

Together, these results demonstrate that random unbiased perturbations of $Y\mid X$ create a new mode of uncertainty: neither adversarial nor per-sample, but dataset-level, and strong enough to shape both practical estimator tuning and minimax rates.

\paragraph{Related work.}
One line of work models distribution shift via distributionally robust optimization (DRO) \citep{blanchet_distributionally_2025, blanchet_optimal_2022, duchi_learning_2021}, treating the deviation between a source $P'$ and a target $P_0$ adversarially, guaranteeing performance in the worst case over all distributions within a divergence ball. This literature assumes a fixed shift and provides robust but pessimistic guarantees. If this method were applied to the RUP setting, the random, mean-zero shifts would be treated as deterministic and biased, which would lead to overly conservative bounds (see \Cref{fig:why_no_dro}). 

A second line of work assumes benign but smooth source-target shifts \citep{cai_transfer_2021, ma_optimally_2023, reeve_adaptive_2021}. This involves a fixed shift with a smooth or bounded density ratio, or a smoothness assumption in the classification boundary. Thus, this literature does not allow for the highly irregular density ratios that arise natural under RUPs, as it controls the shift via regularity rather than randomness.

Closest to our work are recent models that introduce random shifts of the full joint law. These works analyze how random shift models may be used to predict $Y|X$ shift from $X$ shift \cite{jin_beyond_2024}, how to construct confidence intervals that allow for distributional shifts \cite{jeong2022calibrated}, and how to exploit multiple shifted distributions \cite{jeong_out--distribution_2024}. We build on the models introduced in these works by studying the setting setting only allows for $Y|X$ shift and consider the question of optimal convergence rates for nonparametric regression. Unlike classical measurement-error and surrogate-outcome literatures \citep{buhlmann2014high}, which introduce noise at the level of individual observations, our distributional randomness is on the dataset-level as it stems from a single realization $\xi$ per dataset.

Finally, our minimax lower bounds connect to foundational work on nonparametric rates and information-theoretic lower bounds \citep{huber_lower_1997,tsybakov_introduction_2009, yu_assouad_1997}. By isolating random unbiased perturbations of $Y\mid X$, we identify a new mode of uncertainty—neither adversarial nor per-sample—that inflates estimator variance by decreasing the effective sample size from $n$ to $n_\eff$.

The rest of the paper is organized as follows. In \Cref{sec:rups}, we define RUPs and present the partition model (\Cref{ssec:partmodel}) and the correlated noise model (\Cref{ssec:corr_noisemodel}) as two examples for how such RUPs may be generated. In \Cref{sec:ub}, we provide upper bounds for the convergence rate of pointwise risk when sampling from an RUP. To do this we derive \eqref{eq:bv_rup} in the context of local polynomial estimators. In \Cref{sec:lb}, we derive corresponding lower bounds for the correlated noise model. Finally, in \Cref{sec:numres} we provide numerical results from simulated RUPs, which exhibit the expected convergence rates.

\section{Random Unbiased Perturbations (RUPs)}
\label{sec:rups}

We now formalize the notion of \emph{random unbiased perturbations} introduced above. Recall that our motivating examples involved datasets drawn from a single randomly shifted conditional law $P_\xi(Y|X)$ with a fixed covariate distribution $P^X$. In this section, we give a precise definition of such perturbations, specify how their strength and correlation structure are parameterized, and provide two concrete generative examples. These constructions form the foundation for the upper and lower bound results that follow.

\subsection{RUP Definition}
Let $P_0$ be a baseline distribution over $(X,Y)$, where $X \sim P_0^X$, and, setting $f(X) = \E_0[Y|X]$, we write $Y = f(X) + \e$ for some mean zero random $\e$. In this paper, we focus on the setting with the following simplifying assumption (although RUPs could be generalized beyond this).
\begin{assume}
\label{ass:e_indep_x}
    Suppose, under $P_0$, that $X \sim \mathrm{Unif}[0,1]$, and that $Y = f(X) + \e$, where $X \perp \e$, $\E_0[\e]=0$, and $\Var_0(\e) = \sigma^2$. 
\end{assume}
Under RUPs we shift the distribution of $\e|X$, with shifts that may vary across different values of $X$. We let $\rho: \mathcal{X} \times \mathcal{X} \to [-1, 1]$ be a correlation function such that $\rho(x_1, x_2)$ captures how correlated the shifts are at $X = x_1$ and $X= x_2$. Moreover, we define the average correlation parameter
\[
\bar{\rho} := \E_{X_1,X_2}\big[\rho(X_1,X_2)\big],\qquad X_1,X_2\overset{\text{i.i.d.}}{\sim}P_0^X.
\]
\begin{defn}[Random unbiased perturbation of $P_0$]
\label{def:rup-yx}
Let $P_0$ be a baseline law on $(X,Y)$ and let $\Xi$ be a law on perturbation indices $\xi$. For each $\xi$, let $P_\xi$ be a probability measure on $(X,Y)$ (on the same sample space as $P_0$).
We parameterize a shift in $Y|X$ as a shift in $\e := Y - f(X)$ given $X$. Define 
\[
\Delta_\xi(x) := \E_\xi[\e\mid X=x] - \E_0[\e\mid X=x].
\]
We say that the pair $\mathcal R = (\Xi,(P_\xi)_\xi)$ is a \emph{random unbiased perturbation (RUP)} of $P_0$ with variance scale $\delta^2$ and $X$--correlation kernel $\rho$ if the following hold for $P_0^X$–a.e.\ $x,x'$:

\begin{enumerate}[leftmargin=*]
\item \textbf{Fixed marginal.} $P_\xi^X=P_0^X$ a.s.
\item \textbf{Centering.} $\E_\Xi[\Delta_\xi(x)] = 0$
\item \textbf{Variance scale.}
\[
\Var_\Xi(\Delta_\xi(x)) = \delta^2\sigma^2.
\]
\item \textbf{$X$–dependency.}
\[
\Cov_\Xi(\Delta_\xi(x),\,\Delta_\xi(x'))
= \rho(x,x')\,\delta^2\sigma^2.
\]
\end{enumerate}
We say $P_\xi$ is an RUP with strength $\tau = \delta^2 \bar{\rho}$. 
\end{defn}

\begin{remark}[On the correlation kernel]
The function $\rho(x_1,x_2)$ acts as a correlation kernel across covariate values. 
To be consistent with a covariance structure, $\rho$ must be symmetric with $\rho(x,x)=1$, and for every finite set $\{x_i\}_{i=1}^n$, $[\rho(x_i, x_j)]$ is positive semidefinite.
\end{remark}

Finally, we note that first sampling $\xi \sim \Xi$, then $n$ samples from $P_\xi$, is equivalent to sampling from the mixture distribution \begin{equation}
\label{eq:mix_model}
    \bar{P}^n(D)  := \mathbb{E}_\Xi[P_{\xi}^{\otimes n}(D)].
\end{equation}
We now describe an example of how to generate shifts satisfying the RUP definition above.

\subsection{Examples of RUP-generating processes}
\subsubsection{The Partition model}
\label{ssec:partmodel}
\begin{figure}[h]
\label{fig:part_model}
    \centering
    \begin{minipage}[t]{0.42\linewidth}
        \begin{subfigure}[t]{\linewidth}
            \centering
            \includegraphics[width=\linewidth]{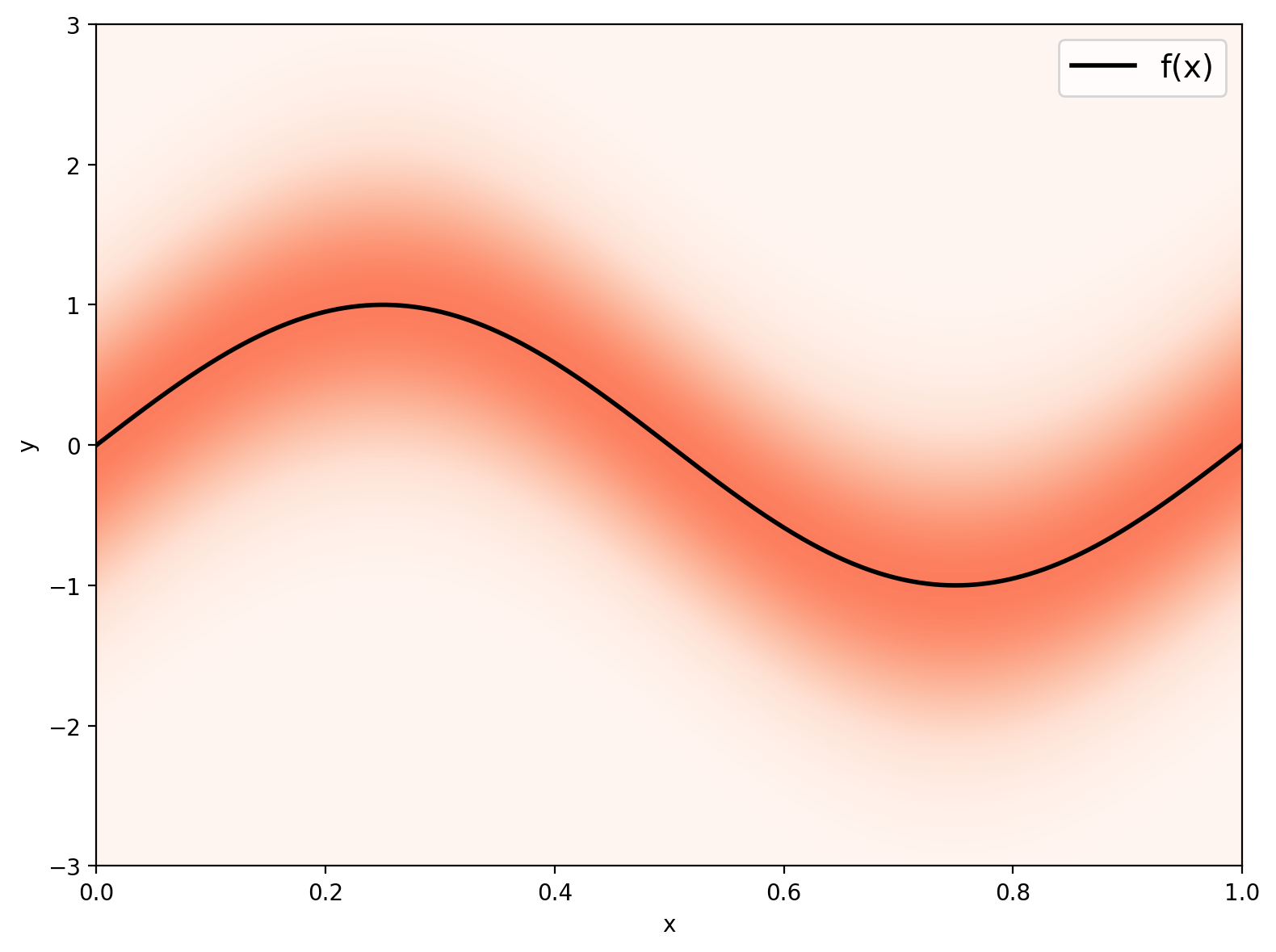}
            \caption{$P_0$: $X \sim \mathrm{Unif}[0,1]$, $Y = f(X) + \mathcal{N}(0, 1/2)$.}
            \label{fig:partmodel:a}
        \end{subfigure}

        \vspace{0.6em}

        \begin{subfigure}[t]{\linewidth}
            \centering
            \includegraphics[width=\linewidth]{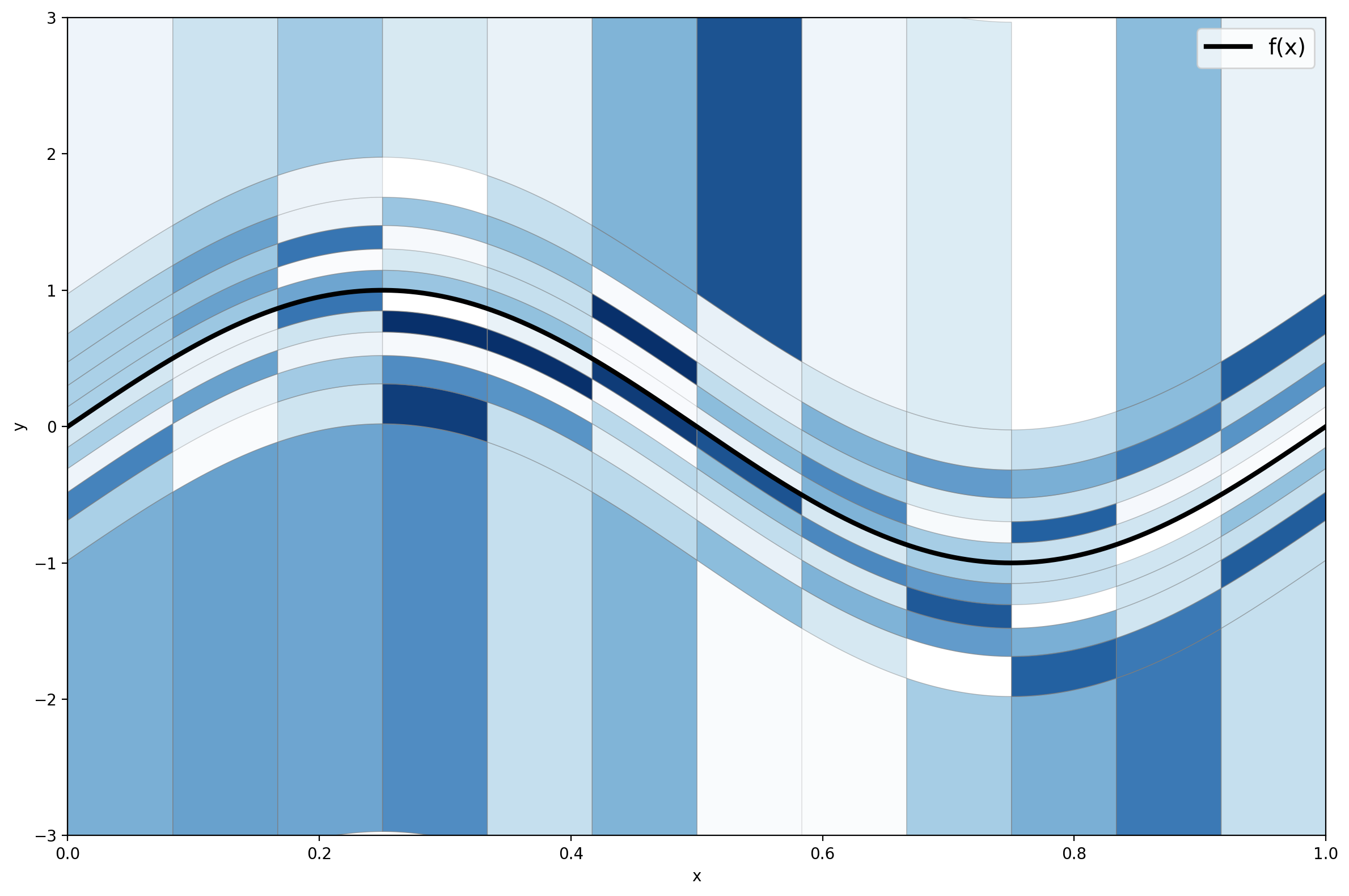}
            \caption{A draw $\xi \sim \Xi$ of random weights, normalized to preserve the marginal $P^X$.}
            \label{fig:partmodel:b}
        \end{subfigure}
    \end{minipage}
    \hfill
    \begin{minipage}[c]{0.52\linewidth}
        \begin{subfigure}[t]{\linewidth}
            \centering
            \includegraphics[width=\linewidth]{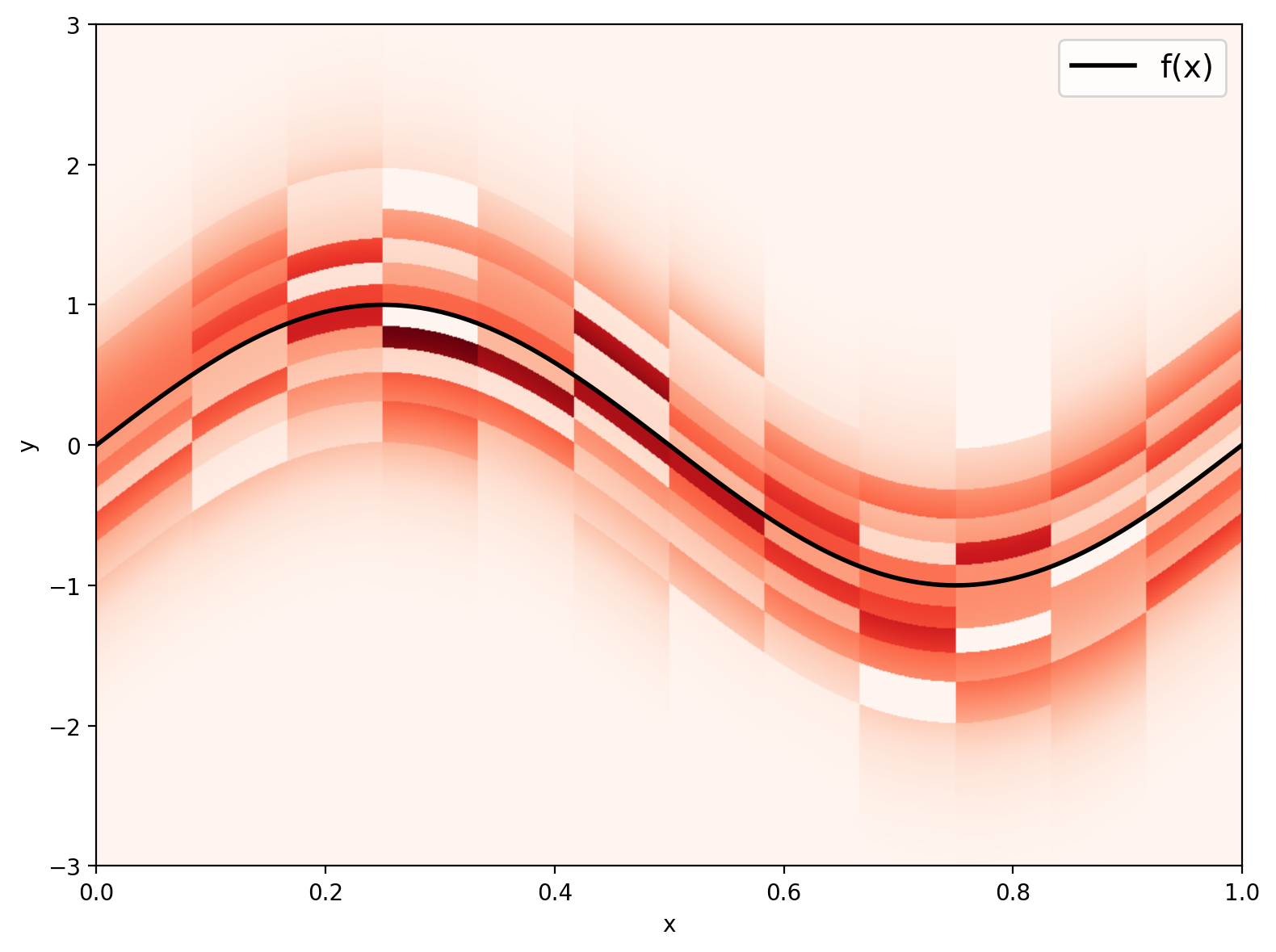}
            \caption{Distribution of $(X,Y)$ under the partition model $P_\xi$: $P_0$ (a) is multiplied by $\xi$ (b).}
            \label{fig:partmodel:c}
        \end{subfigure}
    \end{minipage}
    
    \caption{The partition model. The original density $P_0$ (shown in (a)), is split into $X$ and $\e$ quantiles which are each assigned i.i.d. (normalized) random weights (shown in (b)). The resulting RUP (shown in (c)) is the product of $P_0$ and these weights.}
\end{figure}

We provide an example of an RUP by defining a distribution $\Xi$ over weight arrays $\xi = (\xi_{ij})$, and for each $\xi$ a perturbed law $P_\xi$ on $(X,Y)$. We construct $P_\xi$ by specifying its density in the $(X,\e)$ coordinates and then view it as a law on $(X,Y)$ via $Y=f(X)+\e$. We will not distinguish notationally between the law on $(X,\e)$ and its pushforward to $(X,Y)$.

This model, to which we will refer to as the \textit{partition model}, is based on a construction first introduced by Jeong and Rothenh\"ausler \cite{jeong2022calibrated}. We partition the joint distribution of $(X, \e)$ under $P_0$ into quantiles and generate a random perturbation by randomly upweighting or downweighting the density in each quantile while keeping the marginal distribution over $X$ constant. Formally, we choose integers $B_X,B_\e$,  and set $\{r_i\}_{i=0}^{B_X}$ to be $P_0^X$–quantiles and $\{q_j\}_{j=0}^{B_\e}$ to be $P_0^\e$–quantiles. Then we form bins
\[
   I_i=[r_{i-1},r_i),\quad 
   J_j=[q_{j-1},q_j), \qquad i\in[B_X],\ j\in[B_\e].
\]

Thus under $P_0$, $P_0(X\in I_i)=1/B_X$ and, by $\e\perp X$, $P_0(\e\in J_j\mid X)=1/B_\e$. Next, draw i.i.d.\ positive weights $\{\xi_{ij}\}_{i\le B_X,j\le B_\e}$ (e.g.\ $\xi_{ij}\sim\mathrm{Exp}(1)$).
For $x\in I_i$, $\e \in J_j$ we define the perturbed noise law by
\begin{equation}
    \label{eq:eps-tilt}
    P_{\xi}(x, \e) = \frac{\xi_{ij}}{\frac{1}{\Be}\sum_{j' = 1}^{\Be}\xi_{ij'} }P_0(x, \e)
\end{equation}
In other words, as visualized in \Cref{fig:part_model}, we draw i.i.d. weights for each of $B_X\Be$ quantile bins, then normalize them so that the marginal density of $X$ is unchanged. 
\begin{remark}
Equation~\eqref{eq:eps-tilt} is equivalent to the multiplicative tilt form 
\(
  p_{\xi}(y\mid X=x)=\frac{\xi_{ij}}{\sum_{j'}\xi_{ij'}}\,p_0(y\mid X=x)
\)
whenever $y-f(x)\in J_j$ and $x\in I_i$.
\end{remark}

\begin{prop}
\label{prop:partmodel_valid}
Let $\Xi$ be the law of iid positive weights $\{\xi_{ij}\}_{i \in [B_X], j \in [\Be]}$, and let $\E[\xi^{-3}] < \infty$ (e.g. $\xi \sim \mathrm{Exp}(1)$). The partition model produces a random unbiased perturbation of $P_0$, $(\Xi, (P_{\xi})_\xi)$ with:
\[\text{ Variance scale } \delta^2 = \frac{1}{\Be}\left(\frac{\Var_\Xi (\xi_{ij})}{\E_\Xi [\xi_{ij}]^2} + o(1)\right), \quad \text{X-correlation } \rho(x_1, x_2) = \1_{I_{x_1} = I_{x_2}}, \text{ and }\bar{\rho} = \frac{1}{B_X}.\]
where $I_x$ represents the $P_0^X$-quantile containing $x$.
\end{prop}
See \Cref{sec:bin_rup_proof} for proof.
A dataset $D$ of $n$ samples from a distribution perturbed following the partition model follows the mixture model $\bar{P}^n(D)$ as in \eqref{eq:mix_model}.
\subsubsection{Correlated noise model}
\label{ssec:corr_noisemodel}

\begin{figure}[h]
\label{fig:correlated_noise_model}
    \centering
    \begin{minipage}[t]{0.42\linewidth}
        \begin{subfigure}[t]{\linewidth}
            \centering
            \includegraphics[width=\linewidth]{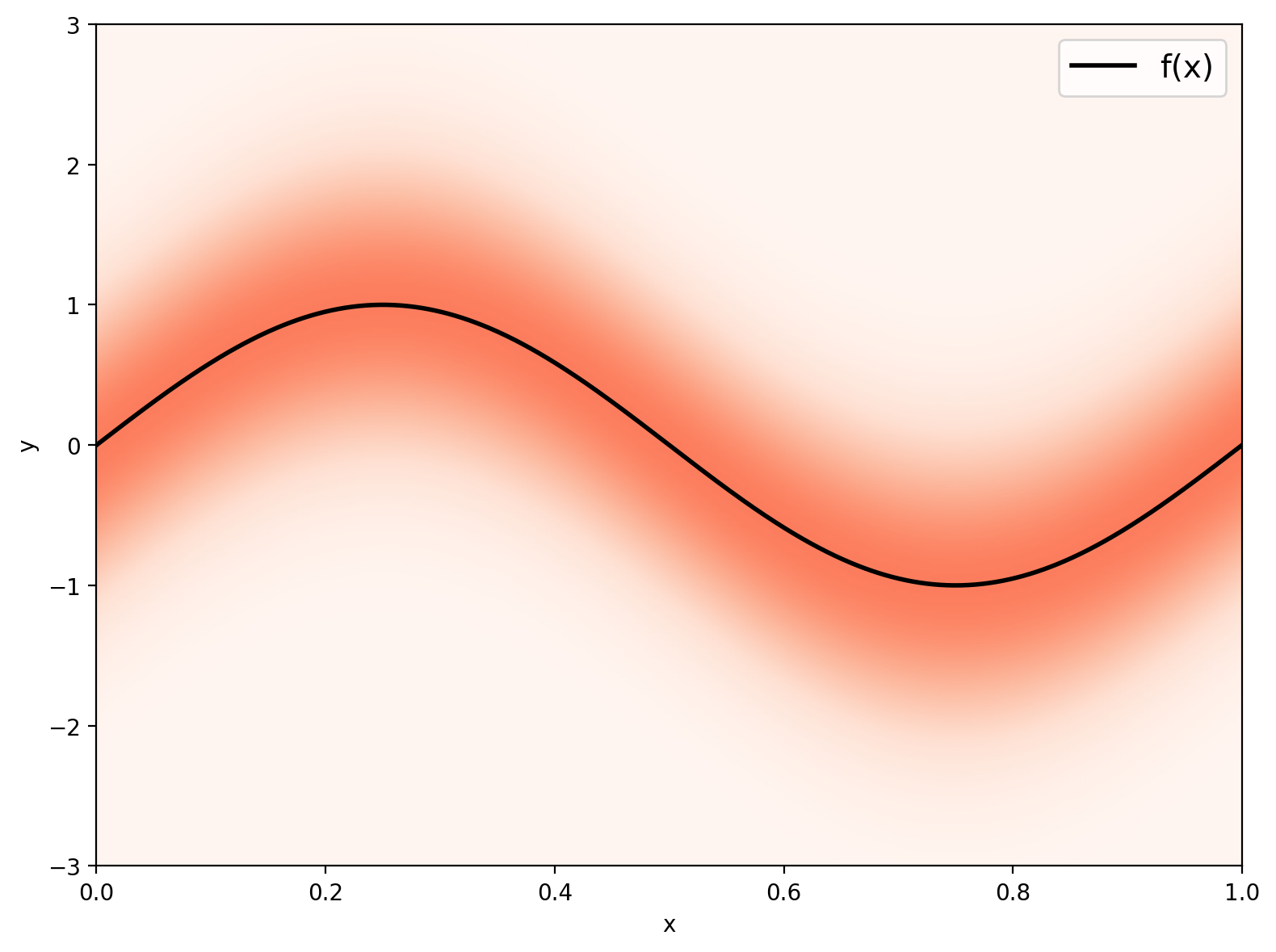}
            \caption{$P_0$: $X \sim \mathrm{Unif}[0,1]$, $Y = f(X) + \mathcal{N}(0, 1/2)$.}
            \label{fig:corrmodel:a}
        \end{subfigure}

        \vspace{0.6em}

        \begin{subfigure}[t]{\linewidth}
            \centering
            \includegraphics[width=\linewidth]{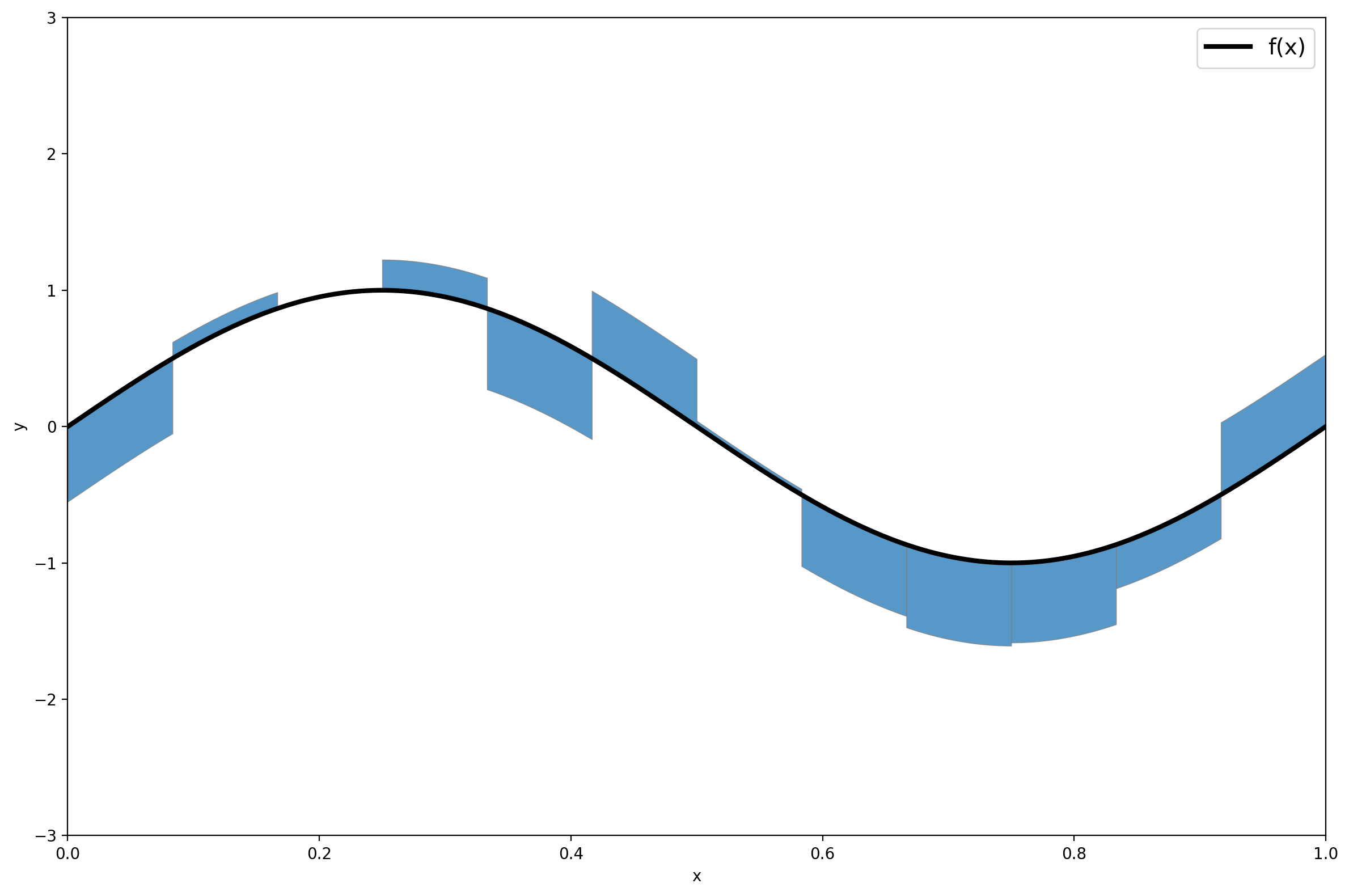}
            \caption{A draw $\xi \sim \Xi$, with $\xi_i \sim \mathcal{N}(0, \delta^2\ \sigma^2)$.}
            \label{fig:corrmodel:b}
        \end{subfigure}
    \end{minipage}
    \hfill
    \begin{minipage}[c]{0.52\linewidth}
        \begin{subfigure}[t]{\linewidth}
            \centering
            \includegraphics[width=\linewidth]{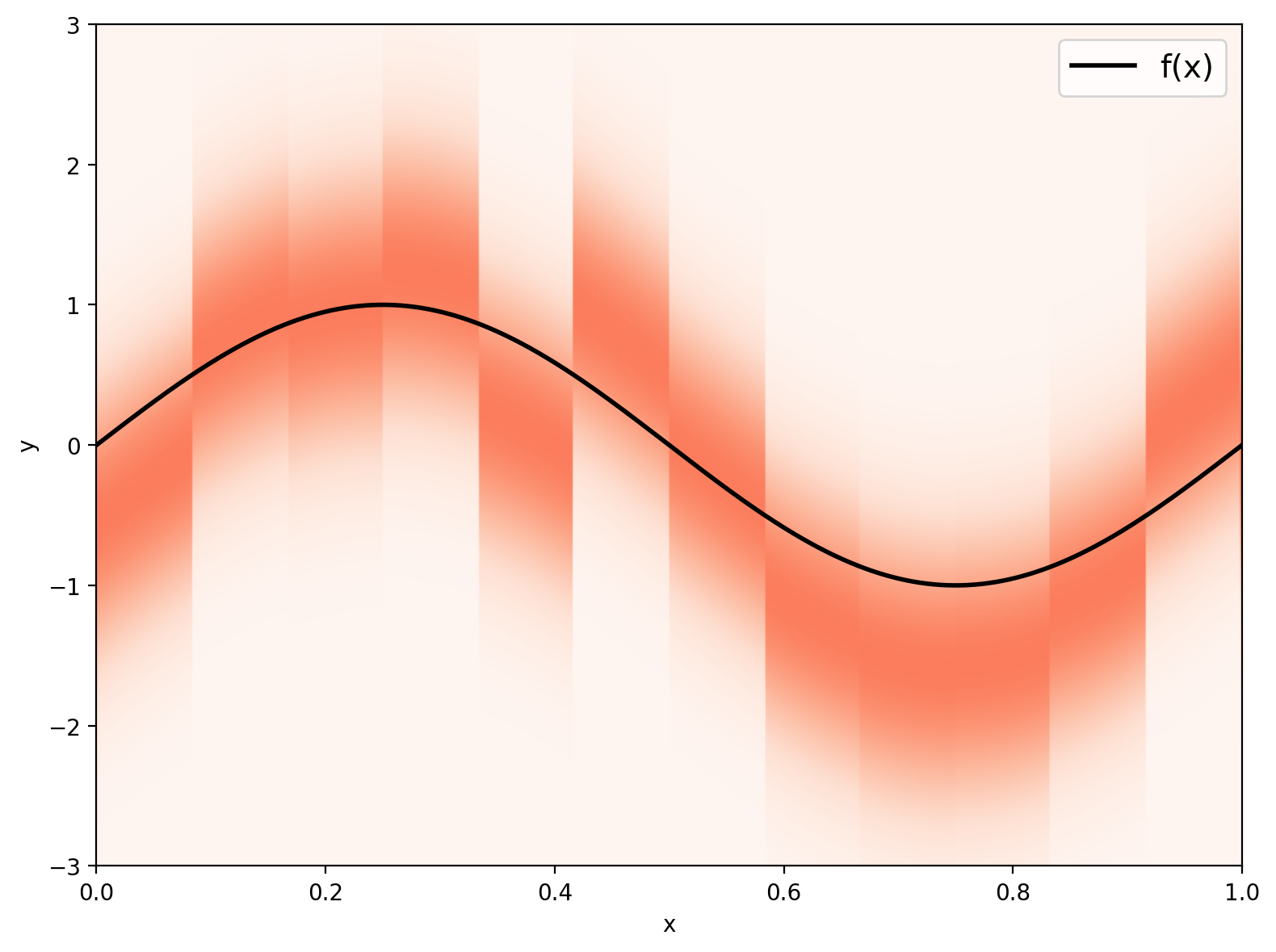}
            \caption{Distribution of $(X,Y)$ under the correlated noise model $P_\xi$: $P_0$ (a) is shifted by $\xi$ (b).}
            \label{fig:corrmodel:c}
        \end{subfigure}
    \end{minipage}
    
    \caption{The correlated noise model. The original density $P_0$ (shown in (a)), is split into $X$-quantiles which are each assigned i.i.d. random shifts (shown in (b)). The resulting RUP (shown in (c)) is the original $P_0$ and shifted by (b).}
\end{figure}

We now introduce an alternative RUP instance, which we will refer to as \textit{the correlated noise model}. Again, we partition the distribution of $X$ into $B_X$ quantiles. This time, instead of multiplicative weights, we have additive noise terms in each quantile.

Recall that $\sigma^2 = \Var_0(\e)$. The perturbation is parameterized by $\xi =\{\xi_b\}_{b=1}^{B_X}$, where $\xi_b \sim^{i.i.d.} \mathcal{N}(0, \sigma^2\delta^2)$.  For observation $X = x$, let $b(x)$ be the $B_X$-bucket and set 
\[P_\xi(x,\e) = P_0^X(x)P_0^\e(\e - \xi_{b(x)})\]
Thus, within every bucket $b$, under $P_\xi$, $Y$ is shifted to be centered at $f(X) + \xi_b$. 
\begin{prop}
\label{prop:corr_noise_valid}
    The correlated noise model produces a random unbiased perturbation of $Y|X$ with:
    \[\text{ Variance scale } \delta^2, \quad \text{X-correlation } \rho(x_1, x_2) = \1_{I_{x_1} = I_{x_2}}, \text{ and }\bar{\rho} = \frac{1}{B_X},\]
where $I_x$ represents the $P_0^X$-quantile containing $x$.
\end{prop}
See \Cref{sec:cor_noise_rup_proof} for proof. 

\section{Upper Bounds under Random Unbiased Perturbations}\label{sec:ub}
In this section, we recall basic properties and convergence rates of local polynomial estimators (LPEs) in the absence of distribution shift. We then analyze how LPEs behave in the RUP setting by deriving the additional distributional variance term. Finally, we discuss the consequences of this additional term for convergence rates and bandwidth tuning. We assume throughout this section that $P_0^X = \mathrm{Unif}[0,1]$.

\subsection{Local polynomial estimator (LPE) Preliminaries}
 Local polynomial (LP) regression offers a flexible and widely used nonparametric method for estimating regression functions. We work under the standard local polynomial design assumptions (LP1)–(LP3) of \cite{tsybakov_introduction_2009}, which ensure that the local Gram matrix is uniformly non-degenerate and the LP weights are well-behaved. Under an i.i.d. design, these conditions hold automatically (up to events of vanishing probability). 
 
 For a given bandwidth $h$, an LP estimator of order $\ell$ uses kernel weighting to locally fit polynomials of degree up to $\ell$, adapting to the structure of the data around each target point. Such estimators admit a linear representation \cite[See ][Chapter 1.6 for details]{tsybakov_introduction_2009} and can be written as 
	\[\hat{f}_n(x) = \sum_{k=1}^n W_k(x) Y_k,\]
	where the weights $W_k(x)$ depend on the target point $x$, covariates $X_1, \dots, X_n$, a kernel function $K$ and the bandwidth $h$. 
In the standard no-shift setting, where samples are drawn i.i.d. from the target distribution, local polynomial estimators achieve minimax-optimal rates for pointwise risk. For a fixed bandwidth $h$, their bias and variance satisfy the following bounds.  

\begin{prop}[based on Tsybakov \cite{tsybakov_introduction_2009}, Proposition~1.13]
\label{prop:bv_noshift}
Let $f \in \Sigma(\beta,L)$ on $[0,1]$ and $\hat f_n(x)$ be the $LP(\ell)$ estimator of $f$. Under standard assumptions (see \Cref{lem:weightprops}), the pointwise risk for any $x_0 \in [0,1]$ admits the following bias variance decomposition
\[
   \E_0\big[(\hat f_n(x_0)-f(x_0))^2\big] = \underbrace{\left(\E_0[\hat{f}_n(x_0)]- f(x_0)\right)^2}_{b^2(x_0;\hat f_n)} + \underbrace{\Var_0(\hat{f}_n(x_0))}_{v^2_{\mathrm{sampling}}(x_0;\hat f_n)},
\]
where
\[
   b^2(x_0;\hat f_n)| \le C_1 h^{2\beta},
   \qquad
   v^2_{\mathrm{sampling}}(x_0;\hat f_n) \le \tfrac{C_2\sigma^2}{n h} + O\left(\frac{h^{2\beta}}{nh}\right),
\]
where $C_1 = \frac{C_* L}{\ell!}$ and $C_2 = \sigma^2 C_*^2$, with $C_*$ depending on the choice of kernel and $\inf_x P_0^X(x)$. The last term in the sampling variance comes from the random design setup. 
\end{prop}
Thus, in the no-shift setting, the pointwise risk is minimized at $h \asymp n^{-1/(2\beta+1)}$, giving the optimal convergence rate 
\[\E_0\big[(\hat f_n(x_0)-f(x_0))^2\big] \le Cn^{-2\beta/(2\beta+1)}\]
for some constant depending on $\ell, L, \sigma^2$ and the choice of kernel.
\subsection{LPEs under RUPs}
We present the following analog to \Cref{prop:bv_noshift}, for when sampling from a RUP.
\begin{theorem}
		\label{thm:ub}
		For any $x_0 \in [0,1]$, let $\hat{f}_n^0(x_0)$, $\hat{f}_n'(x_0)$ be $LP(\ell)$ estimators of bandwidth $h$ using $(X_{1:n}, Y_{1:n})$ from $P_0^{\otimes n}$ and an RUP mixture $\bar{P}^n$ of strength $\tau = \delta^2 \bar{\rho}$ respectively. Moreover, assume that $\rho(X_1, X_2)$ has finite correlation length $\lambda < h$, i.e. $|x_1 - x_2| > \lambda \implies \rho(x_1, x_2) = 0$, and that the choice of $x_0$ is independent of $\rho$. Then
		\begin{equation}
			\E_\Xi\E_{\xi}[(\hat{f}_n(x_0) - f(x_0))^2] \le b^2(x_0; \hat{f}^0_n) + v^2_{\mathrm{sampling}}(x_0;\hat{f}_n^0) + v_{\mathrm{dist}}^2(x_0; \hat{f}'_n),
		\end{equation}
		where
		\[b^2(x_0;\hat{f}_n^0) \le C_1 h^{2\beta}, \quad \quad v^2_{\mathrm{sampling}}(x_0; \hat{f}_n^0) \le \frac{C_2 \sigma^2}{nh}+ O\left(\frac{h^{2\beta}}{nh}\right), \quad \text{ and } v_{\mathrm{dist}}^2(x_0; \hat{f}'_n) \le \frac{C_2 \sigma^2 \tau}{h}+ O\left(\frac{h^{2\beta}}{nh}\right),\]
     where $C_1, C_2$ are as in \Cref{prop:bv_noshift}.
	\end{theorem}
    Note that the distributional variance term $v_{\mathrm{dist}}^2$ has the same bandwidth dependence as the classical $v_{\mathrm{sampling}}^2$ and the new error decomposition is equivalent to replacing the sample size $n$ in \Cref{prop:bv_noshift} with $n_{\eff} =\left(\frac{1}{n} + \tau\right)^{-1} = \frac{n}{1 + n \tau}$

The perturbation strength, $\tau = \delta^2 \bar{\rho}$ quantifies both the variance scale $\delta^2$ of the random perturbations and their spatial correlation $\bar{\rho}$ across covariates. This dependence on $\tau$ contrasts with previous distribution generalization results, which bound the risk in terms of deterministic, pointwise measures of shift. In domain adaptation and DRO frameworks, \cite{ben-david_theory_2010, duchi_learning_2021}, the discrepancy between source and target distributions is fixed and measured globally through $f$-divergences, which, in the case of RUPs will not account for the random cancellations of the model. Consider for instance the KL divergence between $P_0$ and an RUP $P_\xi$ under the partition model.
\[ KL(P_0||P_\xi) = \int p_0(x,\e) \log \frac{p_0(x, \e)}{p_\xi(x, \e)}d(x,\e) = \sum_{b = 1}^{B_X \Be}\log \frac{1}{\xi_b} \int_b p_0(x,\e) d(x,\e) = \frac{1}{B_X \Be}\sum_{b = 1}^{B_X \Be} (- \log \xi_b) \approx \E[-\log \xi]
\]
Thus, as the number of bins $B_X \Be$ increases, this KL divergence converges to the negative log-moment of the weight distribution $\Xi$. However, the perturbation strength $\tau$ tends to zero as $B_X\Be \to \infty$, and indeed as shown by the theorem above, the error tends to that of the no-shift setting.

Similarly, recent transfer-based formulations \cite{reeve_adaptive_2021, sahoo_learning_2025} characterize distribution shift through smooth or bounded transfer functions $T(x) = p_T(x)/p_S(x)$, which capture systematic bias between domains but not stochastic variation of $T(x)$. Our results show that even under irregular changes of the density ratio, transfer is possible.

\subsubsection{Bandwidth selection under RUPs}
\label{sssec:bandwidth_for_rups}
The classical approach to bandwidth selection relies on the i.i.d.\ assumption, where methods such as cross-validation or plug-in estimators balance bias against sampling variance. Under RUPs, however, the additional distributional variance term requires a modified approach.

\paragraph{Oracle bandwidth.} 
When $\tau$ is known, the optimal bandwidth minimizes the total MSE:
\[
h^\star = \arg\min_h \left\{C_1 h^{2\beta} + \frac{C_2}{nh} + \frac{C_2\tau}{h}\right\},
\]
where $C_1, C_2 > 0$ are constants. 
\begin{cor} 
\label{cor:rup_ub_rate}
Balancing the bias term against the combined variance yields \[
h^\star \asymp \left(\frac{1}{n} + \tau\right)^{1/(2\beta+1)} = n_{\mathrm{eff}}^{-1/(2\beta+1)}.
\]
Thus, inserting this bandwidth into the results from \Cref{thm:ub}, we have that the LPE pointwise estimator converges at a rate of $n_{\eff}^{-\frac{2\beta}{2\beta + 1}}$.
\end{cor}
 In the regime where $\tau \gg 1/n$, this simplifies to $\tau^{\frac{2\beta}{2\beta + 1}}$, independent of sample size.

\paragraph{Practical estimation.}
In practice, $\tau$ is unknown and must be estimated from data. A rigorous treatment of optimal estimation procedures for $\tau$ is beyond the scope of this work; here we outline two practical approaches that can be employed depending on the available information.

\emph{1. Domain-structured holdouts.} When the target of inference involves generalization across a known domain structure (e.g., different days, geographical regions, or institutions), practitioners often construct validation sets that respect this structure—holding out entire days, regions, or institutions rather than randomly sampling observations. Our RUP framework provides a model for this common practice: such domain-structured holdout sets naturally capture both sampling and distributional variance, whereas standard random splits only reflect sampling variability. Consequently, bandwidth selection via cross-validation on domain-structured holdouts will appropriately account for the effective sample size $n_{\mathrm{eff}}$ rather than the nominal sample size $n$. This perspective aligns with empirical findings in the domain generalization literature \citep{gulrajani2020search}, where validation procedures that respect the structure of distribution shifts consistently yield better-tuned models.

\emph{2. Estimating $\tau$ from summary statistics.} When summary data from the target distribution are available, one can estimate the distributional uncertainty directly. Suppose we have access to population-level statistics (e.g., means of the outcome) for the target distribution across multiple realizations. By comparing the training data to these summary statistics, we can estimate the strength of the distributional shift $\tau$. Specifically, let $\{ \theta_j\}_{j=1}^J$ be summary statistics from $J$ independent realizations of the target distribution. The variance of $\theta_j$ across realizations provides an estimate of the distributional uncertainty, depending on how the statistic aggregates over $X$. Summaries that preserve global symmetry over $X$ (e.g. averages over values of $X$) vary with the global shift strength $\tau$, while pointwise-in-$x$ summaries vary with $\delta^2$. Fixed, nonvanishing local windows (e.g., an LP equivalent kernel at $x_0$ with bandwidth $h$) target the corresponding local shift factor: $\delta^2$ times the window’s  $X$-correlation—which is exactly the variance inflation entering the LP risk. Thus such summary statistics may be used to calibrate the bandwidth $h^\star$. 

\section{Minimax Lower Bounds under Random Unbiased Perturbations}
\label{sec:lb}
We now show that random unbiased perturbations (RUPs) slow the attainable rate of convergence by effectively reducing the sample size. While the optimal learning rate depends on the exact perturbation mechanism, we show that in the case of the correlated noise model, in the regime where the sampling and distributional uncertainties are of the same order, or the sampling uncertainty dominates, the rate derived in \Cref{cor:rup_ub_rate} is the optimal rate of convergence for pointwise risk. We do so under the following additional assumption.
\begin{assume}
\label{ass:0noise_normal}
    The noise $\e$ is normally distributed under $P_0$, i.e. $P^\e_0 \sim \mathcal{N}(0, \sigma^2)$.
\end{assume}
Recall that under the correlated noise model, a perturbation is parameterized by $\xi = \{\xi_b\}_{b=1}^{B_X}$, where $\xi_b \sim \mathcal{N}(0, \delta^2 \sigma^2)$, and observations follow for $X_i = x$, $Y_i = f(x) + \e$, where $P_{\xi}^\e(\e) = P_0^{\e}(\e - \xi_{b(x)})$. We study the regime where the number of buckets $B_X$ grows with $n$ such that $\lim_{n \to \infty} \frac{n}{B_X(n)} < \infty$.
\begin{remark}
Depending on the smoothness class of functions we assume, the restriction of the regime we study can be weakened. However, due to the discreteness of the RUP examples provided, fixing $B_X$ and sending $n \to \infty$ leads to almost identifiable $P_0$ for sufficiently smooth functions $f$. Thus we recover the behavior of the no-shift setting. If, on the other hand, we fix $n$ and send $B_X \to \infty$, the correlation $\rho \to 0$, and this effectively increases the variance of $\e$, sending the effective sample size back to $n$.
\end{remark}
\begin{theorem}
\label{thm:lb}
Suppose that $\beta > 0$ and $L > 0$. Let $(\Xi, (P_\xi)_{\xi})$ be an RUP following the correlated noise model (defined in \Cref{ssec:corr_noisemodel}) of strength $\tau = \tau_n$. For any fixed point $x_0$, let $\hat{f}_n(x_0)$ be an estimator of $f(x_0)$ constructed using $(X_{1:n}, Y_{1:n})\sim \bar{P}^n$. Let $n_\eff = \frac{n}{1 + n\tau_n}$. Then, as long as $\lim_{n \to \infty} n\tau_n < \infty$, the minimax lower bound satisfies
\begin{equation}
     \liminf_{n \to \infty} \inf_{\hat{f}_n} \sup_{f \in \Sigma(\beta, L)} \mathbb{E}_{\Xi} \mathbb{E}_{\xi, f}^n\left[n_{\eff}^{\frac{2\beta}{2\beta + 1}}\left(\hat{f}_n(x_0) - f(x_0)\right)^2\right] \ge c_1
\end{equation}
for some constant $c_1$ depending on $\sigma^2$, $\beta$ and $L$. 
\end{theorem}
We provide a full proof in \Cref{sec:lb_proof}. The idea follows standard minimax arguments; the key new step is the scaling of the KL divergence under RUPs. In traditional minimax arguments for pointwise risk around $x_0$, one usually constructs two hypothesis functions (see for instance \cite{tsybakov_introduction_2009}), $f_0(x) \equiv 0$ and $f_1(x) = Lh^\beta K\left(\frac{x-x_0}{h}\right)$. Here $K$ is some kernel function and $h$ is the bandwidth. Setting $P_0, P_1$ to be the probability distributions given data is generated via $f_0, f_1$ respectively, one shows that $KL(P_0^{\otimes n}||P_1^{\otimes n}) = C_* n h^{2\beta + 1}$, so letting the bandwidth scale as $h \propto n^{-\frac{1}{2\beta + 1}}$, results in a constant bound on the KL divergence as $n$ grows.

Under RUPs of two close parametric distributions $P_0, P_1$, the usual factorization  breaks because the mixture distribution $\bar{P}_i^n$ no longer factorizes. The next lemma shows how the same expansion holds with $n$ replaced by $n_\eff$. 
\begin{lemma}[KL scaling under the correlated noise model]
\label{lem:KL_neff}
Under the correlated noise model, let $f_0 \equiv 0$, and $f_1(x) = Lh^\beta K\left(\frac{x-x_0}{h}\right)$ for some fixed $x_0$ and fixed kernel function $K$. Let $\bar{P}_0^n, \bar{P}^n_1$ be the mixture probability distributions for the noise-correlated RUP data generated by $f_0$, $f_1$ respectively. Then, there exist a constant $C$, depending on $L$, $\beta$, $K_{\max} = ||K||_\infty$, and $\sigma^2$ such that, in the regime where $\lim_{n \to \infty}\frac{n}{B_X(n)} < \infty$, 
\[\lim_{n \to \infty}\frac{1}{n_\eff}KL(\bar{P}_0^n||\bar{P}_1^n) = Ch^{2\beta +1}\]
In other words, setting the bandwidth to be $h \propto n_{\eff}^{-\frac{1}{2\beta +1}}$ leads to a KL divergence that is bounded by a constant in the limit $n \to \infty$.
\end{lemma}

The lemma is proved in \Cref{subsec:kl_scaling_proof}, but we outline the main argument here. Under the correlated noise RUP, datapoints $(X_i, Y_i)_{i=1}^n$ are no longer iid, however, the $X_i$s are iid, and the $Y|X$ follows a multivariate gaussian distribution, $\mathcal{N}(f(X), \Sigma)$, where 
\[\Sigma_{ij} = \begin{cases}
    (1 + \delta^2)\sigma^2& \text{if }i=j\\
    \delta^2\sigma^2 & \text{if }i \neq j, \text{ but } b(X_i) = b(X_j)\\
    0 & \text{if }b(X_i) \neq b(X_j)
\end{cases}\]
Computing $\Sigma^{-1}$ (see \Cref{subsec:kl_scaling_proof}) shows that the within–bucket average direction is downweighted by a factor $(1 + n\tau)^{-1}$, so the Fisher information grows like an \emph{effective} sample size $n_\eff < n$ rather than $n$ itself. This is exactly the $n_\eff$ appearing in Lemma~\ref{lem:KL_neff}, and explains why the optimal bandwidth scales as $h \asymp n_\eff^{-1/(2\beta+1)}$ in the correlated noise model.
\begin{remark}[Extension to multi-hypothesis minimax arguments.]
The same LAN--based reasoning extends directly to minimax lower bound proofs that rely on more than two hypotheses, such as
Fano, Le~Cam, or Assouad constructions for mean integrated squared error (MISE) risks. 
In those arguments, the pairwise KL divergences for RUPs under the correlated noise model
$\KL(\bar P_i^{(n)}\|\bar P_j^{(n)})$ appearing in the packing or testing steps
scale as $n_{\eff}\,\KL(P_i^0\|P_j^0)$ by \Cref{lem:KL_neff}. 
\end{remark}
\section{Numerical results}
\label{sec:numres}
We illustrate the theoretical results above through simple simulations using an RUP model. For each dataset, covariates $X_i \sim \mathrm{Unif}[0,1]$, outcomes are generated as $Y_i = f(X_i) + \varepsilon_i$ with $f(x) = \sin(2\pi x)$, and $\varepsilon_i$ drawn according to the perturbed law $p_\xi(\varepsilon \mid X)$ following the correlated noise model. We compute the local-polynomial estimator of order~1 with Epanechnikov kernel over a range of bandwidths $h$ and report the empirical mean integrated squared error (MISE) across $100$ perturbation realizations.
\begin{figure}
    \centering
    \includegraphics[width=0.7\linewidth]{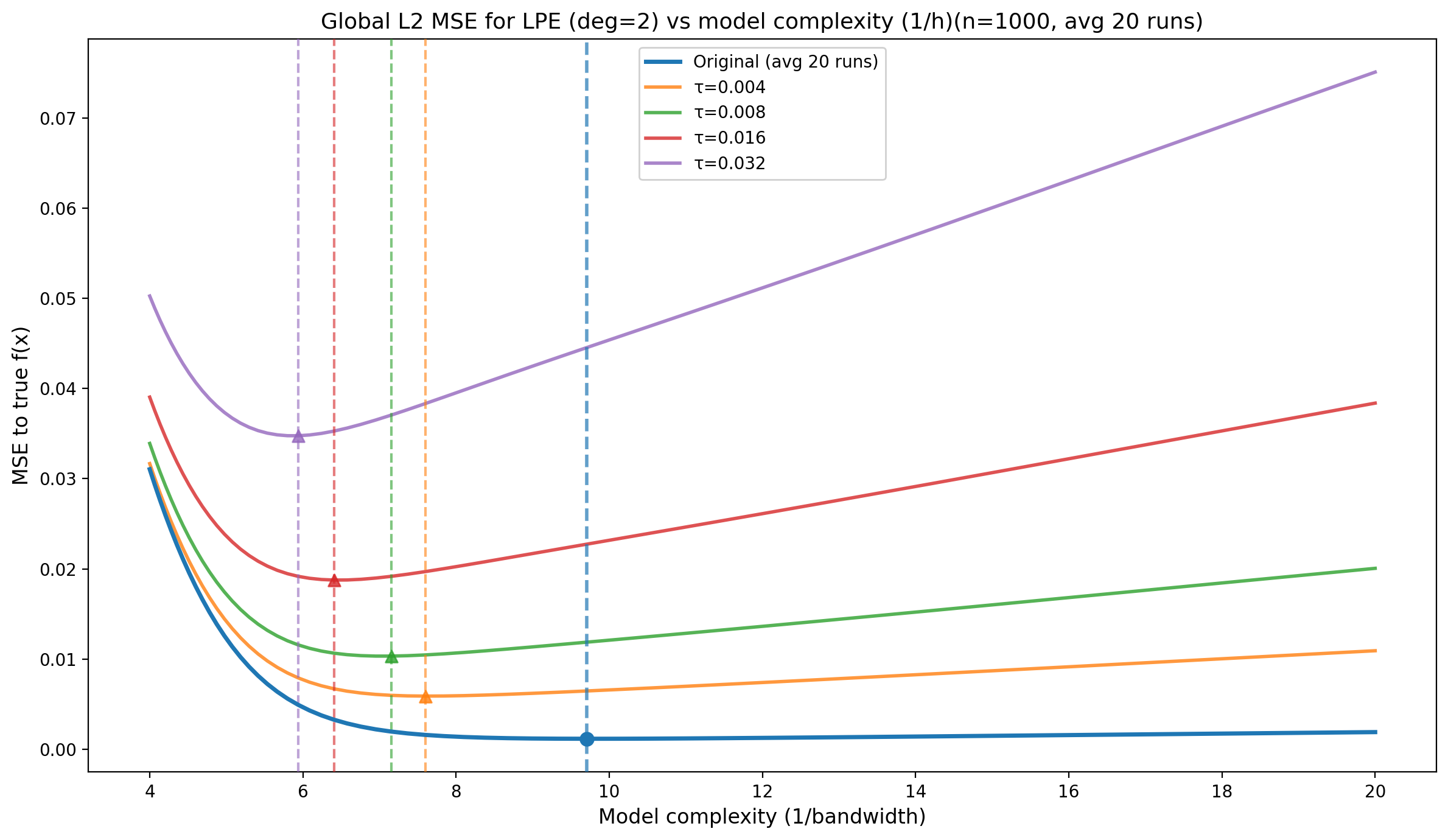}
    \caption{Simulation illustrating the predicted behavior from \Cref{fig:bv_shift}. Under i.i.d. sampling assumptions, a model would choose a smaller bandwidth and hence a higher complexity model. As random distribution shift is introduced, the optimal model complexity decreases to account for the additional variance.}
    \label{fig:mse_w_bandwidth}
\end{figure}

\Cref{fig:mse_w_bandwidth} illustrates the MISE as a function of different bandwidths for three distributions: a baseline $P_0$, and a sequence of perturbed distributions, $\{P_{\xi}\}$. As expected, we see that the bandwidth achieving minimal MISE increases, (i.e. optimal model complexity decreases) as we consider stronger perturbations. In particular, simply tuning the bandwidth size to the shift can lead to significant improvements in the case of larger random shifts.

\begin{figure}
    \centering
    \includegraphics[width=0.7\linewidth]{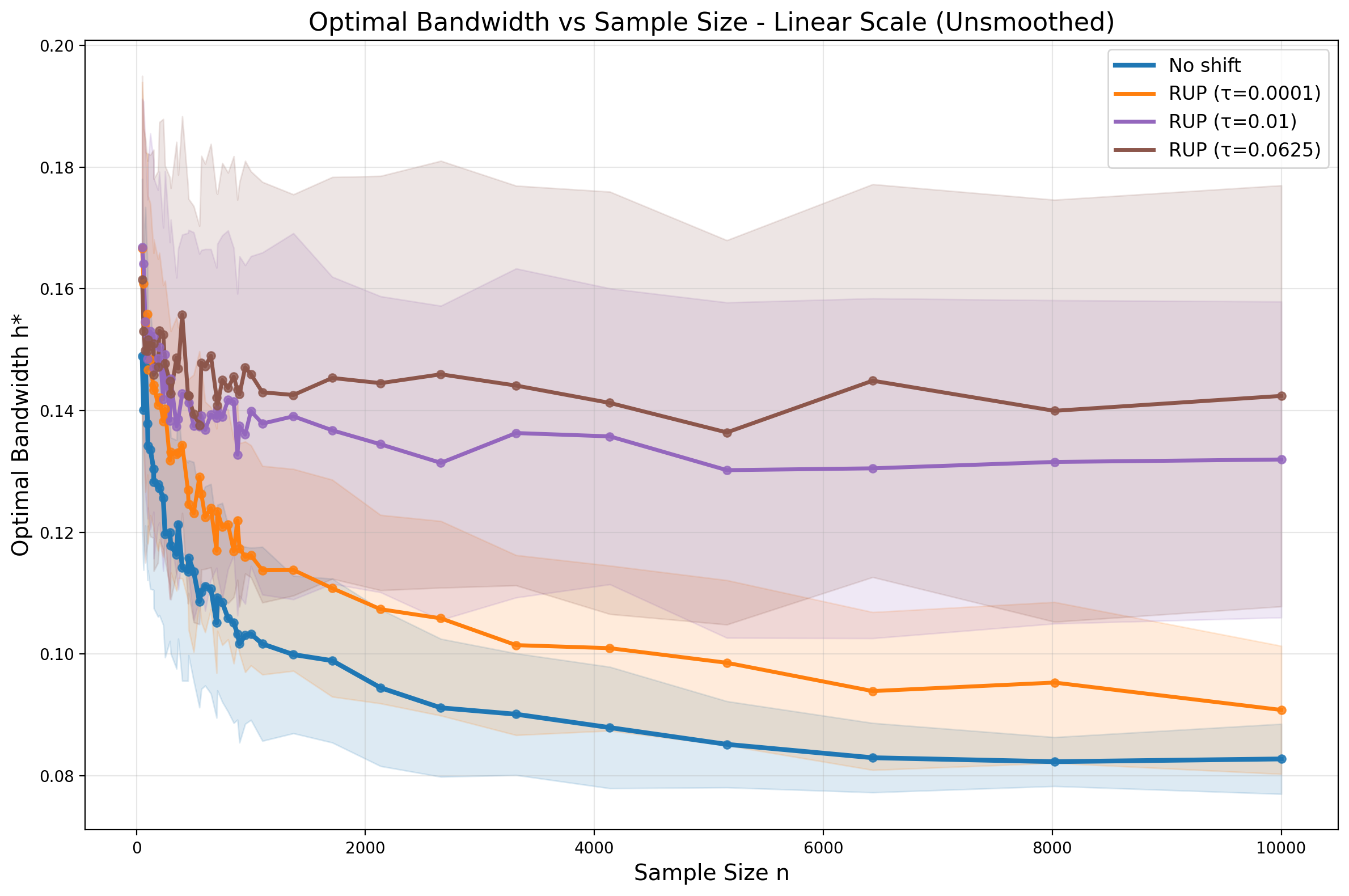}
    \caption{Empirical optimal bandwidth $h^\star$ as a function of sample size $n$ under different perturbation strengths. When $\tau=0$, $h^\star$ decreases as $n^{-1/(2\beta+1)}$ (blue); for $\tau>0$, the curve flattens once $1/n<\tau$, consistent with the theoretical transition to the $\tau^{1/(2\beta+1)}$ scaling.}
    \label{fig:bandwidth_w_n}
\end{figure}

In \Cref{fig:bandwidth_w_n} we analyze how the optimal bandwidth behaves as $n$ grows. Under no shift, the best $h$ is scaled with $n$ and we see it decrease throughout the range of $n$'s considered. However, for perturbed distributions, the optimal bandwidth depends on both $n$ and $\tau$, and once $1/n < \tau$, increasing $n$ no longer leads to smaller optimal bandwidths.

\section{Discussion}
We have introduced a framework for nonparametric regression under random unbiased perturbations (RUPs) of the conditional distribution $Y|X$. This setting captures a common data collection scenario where datasets are drawn from randomly shifted versions of a target distribution, with the marginal distribution of covariates held fixed but the conditional law perturbed by uncontrollable factors such as measurement protocols, environmental conditions, or sampling procedures.

Our analysis reveals that RUPs introduce a distinct mode of uncertainty that differs fundamentally from both classical sampling variability and adversarial distribution shifts. The key insight is that distributional uncertainty manifests as variance inflation rather than bias. This variance inflation reduces the effective sample size to $n_{\mathrm{eff}} = n/(1 + n \tau)$, where $\tau$ quantifies the perturbation strength as the product of variance scale $\delta^2$ and average correlation $\bar{\rho}$.

For local polynomial estimators, we have shown that this effective sample size affects both practical bandwidth selection and fundamental convergence rates. When distributional uncertainty dominates sampling uncertainty ($\tau \gg 1/n$), optimal bandwidths must be chosen to balance bias against distributional variance rather than sampling variance, leading to the scaling $h \propto \tau^{1/(2\beta+1)}$ and convergence rate $\tau^{2\beta/(2\beta+1)}$. Our minimax lower bounds demonstrate that these rates are optimal, at least in the setting $\lim_{n \to \infty} n \tau_n < \infty$, by studying the correlated noise model and establishing that distributional uncertainty fundamentally limits the rate of convergence.

The RUP framework offers several advantages over existing approaches to distribution shift. Unlike DRO or transfer-learning approaches, which assume fixed or adversarial shifts, our framework models dataset-level randomness that averages to the target law, yielding sharper and more realistic risk characterizations. Unlike transfer learning and domain adaptation settings, which treat the source-to-target shift as fixed, we model dataset-level randomness that averages to the target distribution. This perspective is more appropriate for applications where perturbations represent idiosyncratic variation rather than systematic domain differences.

Several directions for future work emerge from our analysis. First, while we have focused on local polynomial estimators, the effective sample size principle should extend to other nonparametric methods such as kernel ridge regression, smoothing splines, and nearest-neighbor estimators. From a methodology perspective, it would also be interesting to see whether the sample size principle can be applied to random forest, gradient boosting, or deep learning techniques. Second, our bandwidth selection discussion has emphasized the theoretical optimum; developing robust procedures that estimate $\tau$ based on principles discussed in \Cref{sssec:bandwidth_for_rups} remains an important part of the puzzle.
From the lower bound side, we have focused on perturbation models that are discrete in that they depend on quantile binning of the distribution and therefore only allow lower bounds in the regime where $\lim_{n \to \infty} n\tau_n < \infty$. However, there may be continuous RUPs based on random processes that could extend these bounds to all regimes.

Moreover, extending the framework to structured perturbations with known correlation patterns could broaden applicability. More generally, one could decompose real-world distribution shifts into systematic and random components, thus capturing both persistent biases \textbf{across domains} and transient, mean-zero fluctuations \textbf{within domains}. Finally, similar variance inflation phenomena may occur in other statistical problems, such as density estimation, classification, or high-dimensional inference. Investigating these would deepen our understanding of how dataset-level randomness affects statistical inference more broadly.

In summary, random unbiased perturbations provide a tractable model that captures plausible patterns of distributional uncertainty and leads to concrete modifications of classical nonparametric methods. By recognizing that such uncertainty inflates variance rather than bias, practitioners can make more informed choices about model complexity and achieve better performance than methods tuned under idealized i.i.d.\ assumptions or overly conservative distributional robustness constraints.

\section{Acknowledgments}

Rothenh\"ausler gratefully acknowledges support as a David Huntington Faculty Scholar, Chamber Fellow, and from the Dieter Schwarz Foundation.

\pagebreak

\bibliography{references}
\bibliographystyle{abbrv}
\pagebreak
\appendix
\section{Validity of RUP models}
\subsection{Validity of the partition model}
\label{sec:bin_rup_proof}
Here we prove \Cref{prop:partmodel_valid}, that the partition model does indeed define a family of RUPs satisfying \Cref{def:rup-yx}.
\begin{proof}[Proof of \Cref{prop:partmodel_valid}]
 Recall that for $x\in I_i$, $\e \in J_j$ we define the perturbed noise law by
 \[P_\xi(x, \e) = \frac{\xi_{ij}}{\bar{\xi}_i}P_0(x, \e),\]
 where $\bar{\xi}_i = \frac{1}{\Be}\sum_{j' = 1}^{\Be}\xi_{ij'}$.
\begin{enumerate}[leftmargin=*]
\item \textbf{Fixed marginal.} We write out for any $x$, wlog assume $x \in I_i$,
    \begin{align*}
        P_\xi(x) &= \int P_\xi(x, \e)d\e = \sum_{j=1}^{\Be} \int_{J_j} P_\xi(x, \e)d\e\\
        &= \sum_{j=1}^{\Be} \int_{J_j} \frac{\xi_{ij}}{\bar{\xi}_i }P_0(x, \e)d\e
        \intertext{By \Cref{ass:e_indep_x},}
        &= \sum_{j=1}^{\Be}\frac{\xi_{ij}}{\bar{\xi}_i } \int_{J_j} P_0(x) P_0(\e)d\e\\
        &= \sum_{j=1}^{\Be}\frac{\xi_{ij}}{\bar{\xi}_i } \frac{1}{\Be} P_0(x) = P_0(x)
    \end{align*}
    \item \textbf{Centering.} Recall $\Delta_\xi(x) := \E_\xi[\e\mid X=x] - \E_0[\e\mid X=x]$. Then, under the partition model, 
    \begin{align*}
        \Delta_\xi(x) &= \int \e \left(P_\xi(\e|X = x) - P_0(\e |X =x)\right)d\e\\
        &= \sum_{j=1}^{\Be} \int_{J_j} \e \left(\frac{\xi_{ij}}{\bar{\xi}_i } - 1 \right) P_0(\e |X =x)d\e\\
        &= \frac{1}{\Be}\sum_{j=1}^{\Be} \left(\frac{\xi_{ij}}{\bar{\xi}_i } - 1 \right) \E_0[\e|\e \in J_j]
    \end{align*}
    Now, taking expectations over $\Xi$, we have $\E_\Xi\left[\frac{\xi_{ij}}{\bar{\xi}_i } - 1 \right] = 0$ for all $i \in [B_X], j \in [\Be]$.
\item \textbf{Variance scale.} Let $m_j = \E_0[\e|\e \in J_j]$. From above, we have
\begin{align*}
    \Var_\Xi\left(\Delta_\xi(x)\right) &= \frac{1}{\Be^2}\sum_{j=1}^{\Be} \Var_\Xi \left(\frac{\xi_{ij}}{\bar{\xi}_i }\right) m_j^2 + \frac{1}{\Be^2}\sum_{\substack{1 \le j,k\le \Be \\ j \neq k}}\Cov_\Xi \left(\frac{\xi_{ij}}{\bar{\xi}_i }, \frac{\xi_{ik}}{\bar{\xi}_i }\right) m_j m_k\\
    \intertext{Now by Taylor expansion since $\E[\bar{\xi}^{-3}] \le \E[\xi^{-3}]< \infty$, we can write $\Var_\Xi \left(\frac{\xi_{ij}}{\bar{\xi}_i }\right) = \frac{\Var(\xi)}{\E[\xi]^2} (1 - O(1/\Be))$. Moreover, for $j \neq k$, $\Cov_\Xi \left(\frac{\xi_{ij}}{\bar{\xi}_i }, \frac{\xi_{ik}}{\bar{\xi}_i }\right) = - O(1/\Be)\frac{\Var(\xi)}{\E[\xi]^2} $. Thus}
    \Var_\Xi\left(\Delta_\xi(x)\right) &= \frac{\Var(\xi)}{\E[\xi]^2} \frac{1}{\Be^2}\left(\sum_{j=1}^{\Be} m_j^2 (1 - O(1/\Be)) - O(1/\Be) \left(  \left(\sum_{1 \le j\le \Be} m_j \right)^2 - \sum_{j=1}^{\Be} m_j^2\right) \right)\\
    \intertext{Note that $\sum_{1 \le j\le \Be} m_j = 0$}
    &= \frac{\Var(\xi)}{\E[\xi]^2} \frac{1}{\Be^2}\left(\sum_{j=1}^{\Be} m_j^2 (1 - O(1/\Be)) \right)\\
\end{align*}
Finally, we will show that $\frac{1}{\Be}\sum_{j=1}^{\Be} m_j^2  = \Var_0(\e) + o(1)$, which will give us the desired result.  
\begin{equation*}
    \frac{1}{\Be}\sum_{j=1}^{\Be} m_j^2 = \E_0[\E_0[\e|J]^2] = \Var_0(\E_0[\e|J])+\underbrace{\E_0[\E_0[\e|J]]^2}_{0} = \Var_0(\e) -\underbrace{\E_0[\Var_0(\e|J)]}_{R_{\Be}},
\end{equation*}
where since $\E[\e|J] \to \e$ in $L^2$, we also have that $R_{\Be} = \E_0[(\e - \E_0[\e|J])^2] \to 0$ as $\Be \to \infty$. 
\item \textbf{X-dependency.} We write out the covariance. Let $x_1 \in I_{i_1}$ and $x_2 \in I_{i_2}$
\begin{align*}
    \Cov_\Xi \left(\Delta_\xi(x_1), \Delta_\xi(x_2)\right) &= \Cov_\Xi \left(\frac{1}{\Be}\sum_{j=1}^{\Be} \left(\frac{\xi_{i_1 j}}{\bar{\xi}_{i_1}} - 1\right)m_j, \frac{1}{\Be}\sum_{j=1}^{\Be} \left(\frac{\xi_{i_2 j}}{\bar{\xi}_{i_2}} - 1\right)m_j\right)\\
    &= \begin{cases}
        \Var_\Xi(\Delta_\xi(x_1)) & \text{if }i_1 =i_2\\
        0 & \text{ otherwise}
    \end{cases}
\end{align*}
By previous parts, we are done.
\end{enumerate}
\end{proof}

\subsection{Validity of the correlated noise model}
\label{sec:cor_noise_rup_proof}
Here we prove \Cref{prop:corr_noise_valid}, that the correlated noise model does indeed define a family of RUPs satisfying \Cref{def:rup-yx}.

\begin{proof}[Proof of \Cref{prop:corr_noise_valid}]
 Recall that under the correlated noise model, we observe
 \[P_\xi(x,\e) = P_0^X(x)P_0^\e(\e - \xi_{b(x)})\]
 with $\xi_{b(x)} \sim \mathcal{N}(0, \delta^2 \sigma^2)$, where $\sigma^2 = \Var_0(\e)$.
\begin{enumerate}[leftmargin=*]
\item \textbf{Fixed marginal.} This holds by definition. For $x \in I_i$,
\[P_\xi(x) = \int P_\xi(x, \e)d\e = \int P_0^X(x)P_0^\e(\e - \xi_{b(x)})d\e = P_0^X(x)\]
    \item \textbf{Centering.} Recall $\Delta_\xi(x) := \E_\xi[\e\mid X=x] - \E_0[\e\mid X=x]$. Then,
    \[\Delta_\xi(x) = \int \e' dP_0^\e(\e' - \xi_{b(x)}) - \int\e dP_0^\e(\e) = \int (\e + \xi_{b(x)} -\e)dP_0^\e(\e)  = \xi_{b(x)}\]
    Now, taking expectations over $\Xi$, we have for all $x$, $\E_\Xi\left[\xi_{b(x)}\right] = 0$.
\item \textbf{Variance scale.} From above, we have
\[\Var_\Xi\left(\Delta_\xi(x)\right) = \Var_\Xi\left(\xi_b\right) = \delta^2 \sigma^2\]
\item \textbf{X-dependency.} We write out the covariance. Let $x_1 \in I_{i_1}$ and $x_2 \in I_{i_2}$
\[\Cov_\Xi \left(\Delta_\xi(x_1), \Delta_\xi(x_2)\right)= \Cov_\Xi(\xi_{b_1}, \xi_{b_2})= \mathbf{1}_{I_{i_1} = I_{i_2}} \delta^2 \sigma^2\]
\end{enumerate}
\end{proof}

\section{Proof of Upper Bounds (\Cref{thm:ub})}

We outline the following properties of the weights, which we will use in our analysis. These follow under mild regularity conditions, including a minimal eigenvalue condition on the local design matrix as in \cite[Section 1.6.1]{tsybakov_introduction_2009}.
	
	\begin{lemma}[Tsybakov \cite{tsybakov_introduction_2009}, Lemma 1.3]
		\label{lem:weightprops}
		Assume the kernel function $K$ is bounded by some $K_{\max}$ and has compact support in $[-1,1]$. Moreover, assume that the local design matrix is uniformly non-degenerate (via some $\lambda_0$) for large enough $n$, and the points $X_i$ are dense in $[0,1]$. Then, almost surely,
		\begin{enumerate}[label=(\roman*)]
			\item $\sum_{k=1}^n W_k(x) =1$, with $\sum_{k=1}^n |W_k(x)| \le C_*$
			\item $\sup_{k,x}|W_k(x)| \le \frac{C_*}{nh}$
			\item $W_k(x) = 0$ if $|X_k - x| > h$
		\end{enumerate}
		where $C_*$ is a constant depending only on $\lambda_0$ and $K_{\max}$.
	\end{lemma}

    The assumptions required for the lemma above hold almost surely for the setting with iid distributed $X_i$, where the probability density of $X$ over $[0,1]$ is bounded away from 0. The minimum value of this density will determine $\lambda_0$. Moreover, in this setting we will satisfy the assumption that the points $X_i$ are dense in $[0,1]$ almost surely.

\begin{remark}[Regarding \Cref{prop:bv_noshift}]
    While this proposition is based on Proposition 1.13 from \cite{tsybakov_introduction_2009}, the version in the book only considers the fixed design regime. To extend the result to random design regime, we must account for the variance from $X$ which is the same expression as we consider in \eqref{eq:random_design_variance} and is the source of the $O\left(\frac{h^{2\beta}}{nh}\right)$ term.
\end{remark}
\begin{proof} [Proof of \Cref{thm:ub}] 
By the same decomposition as in \eqref{eq:bv_rup}, we decompose the risk into bias, sampling variance and distributional variance by the law of total variance:
\begin{align*}
    R_n(x_0; \hat{f}) &= \E_\Xi\E_{\xi}\left[(\hat f_n^\xi(x_0)-f^0(x_0))^2\right]\\
    &= \underbrace{\left(\E_\Xi\E_{\xi}\left[\hat f_n (x_0)\right] - f(x_0)\right)^2}_{\text{bias squared}} + \underbrace{\E_\Xi\Var_{\xi}\left(\hat f_n^\xi(x_0)\right)}_{\text{sampling variance}} + \underbrace{\Var_\Xi\E_{\xi}\left[\hat f_n^\xi(x_0)\right]}_{\text{distributional variance}}
    \end{align*}
    The first two terms follow standard LPE analysis under random design setting. Recall that $\hat{f}_n (x_0)= \sum_{k=1}^n W_k(x_0)Y_k$. Then, for the \textbf{squared bias}, we use the centering property and standard LPE analysis (see \cite{fan_local_1996, tsybakov_introduction_2009}).
    \begin{align*}
    \left(\E_\Xi\E_{\xi}\left[\hat f_n (x_0)\right] - f(x_0)\right)^2 &= \left(\E_\Xi\E_{\xi}\left[\sum_{k=1}^n Y_k W_k(x_0)\right] - f(x_0)\right)^2\\
    &= \left(\E_0\left[\sum_{k=1}^n W_k(x_0)f(X_k)\right] - f(x_0)\right)^2 = \frac{C_*^2 L^2}{(\ell!)^2}h^{2\beta}
\end{align*}

For the \textbf{sampling variance}, we have
\[
    \E_\Xi\Var_{\xi}\left(\hat f_n^\xi(x_0)\right) = \underbrace{\E_\Xi \E_X \left[\Var_{\xi}\left(\hat f_n^\xi(x_0) |X \right)\right]}_{I^s} + \underbrace{\E_\Xi \Var_X \left(\E_{\xi}\left[\hat f_n^\xi(x_0) |X \right]\right)}_{II^s}\]
    Note that given $\xi, X$, $\e_k$ are independent under $\mathbb{P}_\xi$.
    \begin{align*}
    I^s &= \E_\Xi \E_X \left[\Var_{\xi}\left(\sum_{k=1}^n W_k(x_0) Y_k |X \right)\right] = \E_\Xi \E_X \left[\Var_{\xi}\left(\sum_{k=1}^n W_k(x_0) \e_k |X \right)\right]\\
    &= \E_\Xi \E_X \left[\sum_{k=1}^n W_k(x_0)^2  \Var_{\xi}\left(\e_k |X \right)\right] = \E_X \left[\sum_{k=1}^n W_k(x_0)^2  \E_\Xi \left[\E_{\xi}\left[\e_k^2 |X \right] - \E_{\xi}\left[\e_k |X \right]^2\right]\right]\\
    &= \E_X \left[\sum_{k=1}^n W_k(x_0)^2 \left( \E_0\left[\e_k^2 |X \right]  - \Var_\Xi (\E_{\xi}\left[\e_k |X \right])\right)\right]\\
    &\le \sum_{k=1}^n\E_X \left[ W_k(x_0)^2 \right]\sigma^2 \le nh\left(\frac{C_*}{nh}\right)^2 \sigma^2 = \frac{\sigma^2 C_*^2}{nh},
\end{align*}
where we use \Cref{lem:weightprops} and $C_*$ depends only on $\lambda_0$ and $K_{\max}$.
For the second, smaller order term (which comes from the choice to include a random design setup), we have 
\begin{equation}
\label{eq:random_design_variance}
    II^s = \E_\Xi \Var_X \left(\E_{\xi}\left[\hat f_n^\xi(x_0) |X \right]\right) = \E_\Xi \Var_X \left(\sum_{k=1}^n W_k(x_0) f(X_k)\right) = O\left(\frac{h^{2\beta}}{nh}\right)
\end{equation} 
Finally, for the \textbf{distributional variance},  
\[Var_\Xi \E_{\xi}\left[\hat f_n^\xi(x_0)\right]  = \underbrace{\E_X \Var_\Xi \left(\E_{\xi}\left[\hat f_n^\xi(x_0)|X\right]\right)}_{I^D} + \underbrace{\Var_X \left(\E_\Xi \E_{\xi}\left[\hat f_n^\xi(x_0)|X\right]\right)}_{II^D}\]
For the first, main term, we have:
\begin{align*}
    I^D &= \E_X \Var_\Xi \left(\E_{\xi}\left[\hat f_n^\xi(x_0)|X\right]\right)\\
    &= \E_X \Var_\Xi \left(\sum_{k=1}^n W_k(x_0)\E_\xi[\e_k|X_k]\right)\\
    &= \E_X \left[ \sum_{k, \ell} W_k(x_0) W_\ell(x_0) \Cov_\Xi\left(\E_\xi[\e_k|X_k], \E_\xi[\e_\ell|X_\ell] \right) \right]\\
    &= \delta^2 \sigma^2  \E_X \left[ \sum_{k, \ell} W_k(x_0) W_\ell(x_0) \rho(X_k, X_\ell) \right]\\
    &= \delta^2 \sigma^2  \E_X \left[ \sum_{k: |X_k - x_0| < h} W_k(x_0) \sum_{\ell : |X_k - X_\ell| < \lambda}W_\ell(x_0) \rho(X_k, X_\ell) \right]
    \end{align*}
    Since we sum over $W_k(x_0)$ where $|X_k - x_0| <h$, we are summing about $nh$ terms, and for each of these, we sum over $\ell$ such that $|X_k - X_\ell| < \lambda$ which is at most about $n\lambda$ terms. Moreover, for all $k$, $W_k(x)| \le \frac{C_*}{nh}$. Hence, (using the fact that the choice of $x_0$ is independent from the correlation kernel $\rho$),
    \begin{align*}
    &\le \delta^2 \sigma^2  nh \cdot n\lambda \left(\frac{C_*}{nh}\right)^2\E_X \left[\rho(X_k, X_\ell) \mid |X_k - X_\ell| < \lambda\right]\\
    &\le \frac{\delta^2 \sigma^2 C_*^2 \bar{\rho}}{h} = \frac{C_*^2 \sigma^2 \tau}{h} ,
\end{align*}
where the last line follows from
\begin{align*}
    \bar{\rho} &= \mathbb{P}(|X_1 - X_2| \ge \lambda])\underbrace{\E_X[\rho(X_1, X_2)| |X_1 - X_2| \ge \lambda]}_{0} + \mathbb{P}(|X_1 - X_2|< \lambda])\E_X[\rho(X_1, X_2)| |X_1 - X_2|< \lambda]\\ 
    &\ge \lambda\E_X[\rho(X_1, X_2)| |X_1 - X_2|< \lambda]
\end{align*}
Again, for the second, smaller order term, 
\begin{align*}
    II^D &= \Var_X \left(\E_\Xi \E_{\xi}\left[\hat f_n^\xi(x_0)|X\right]\right) = \Var_X \left(\sum_{k=1}^n W_k(x_0)(f(X_k) + \E_\Xi \E_{\xi}\left[\e|X\right])\right) = \Var_X \left(\sum_{k=1}^n W_k(x_0)f(X_k)\right) = O\left(\frac{h^{2\beta}}{nh}\right)
\end{align*}
\end{proof}

\section{Proof of Lower Bounds (\Cref{thm:lb})}
\label{sec:lb_proof}

By Markov's inequality, we write out 
\begin{align*}
    \inf_{\hat{f}_n} \sup_{f \in \Sigma(\beta, L)}\E_\Xi \E_{\xi, f, n}[|\hat{f}_n(x_0) - f(x_0)|] &\ge \inf_{\hat{f}_n} \sup_{f \in \Sigma(\beta, L)}s \E_{\Xi}\mathbb{P}_{\xi, f}^n(|\hat{f}_n(x_0) - f(x_0)| \ge s)\\
     &\ge s \inf_{\hat{f}_n} \max_{f \in \{f_0, f_1\}}\bar{\mathbb{P}}_{\xi, f}^n(|\hat{f}_n(x_0) - f(x_0)| \ge s)
\end{align*}
For any two $f_0, f_1 \in \Sigma(\beta,L)$. In particular, if we set $f_0(x) \equiv 0$, and $f_1(x) = Lh^\beta K\left(\frac{x - x_0}{h}\right)$, where we choose $K$ to be a bounded kernel function satisfying $K \in \Sigma(1/2, L) \cap C^{\infty}$, and vanishing outside of $[-1/2, 1/2]$. Let $\psi(D) = \arg\min_{i \in \{0,1\}}|\hat{f}_n(D, x_0) - f_i(x_0)|$ be the minimal distance test, i.e. returns the index of the function in $\{f_0, f_1\}$ closest to the estimator. Note that in this case, $|f_0(x_0) - f_1(x_0)| = Lh^\beta K_{max}$, so setting $s = \frac{Lh^\beta}{2}K_{\max}$ guarantees that the event $\{\psi(D) \neq j\} \subseteq \{|\hat{f}_n(x_0) - f(x_0)| >s\}$ because there exists a $k$ such that $|\hat{f}_n(x_0) - f_k(x_0)|\le |\hat{f}_n(x_0) - f_j(x_0)|$, thus
\[2s \le |f_k(x_0) - f_j(x_0)| \le |\hat{f}_n(x_0) - f_k(x_0)| + |\hat{f}_n(x_0) - f_j(x_0)| \le 2|\hat{f}_n(x_0) - f_j(x_0)|\]
Thus we can write
\[\inf_{\hat{f}_n} \sup_{f \in \Sigma(\beta, L)}\E_\Xi \E_{\xi, f, n}[|\hat{f}_n(x_0) - f(x_0)|] \ge \frac{Lh^\beta}{2} \underbrace{\inf_{\hat{f}_n} \max_{f_j \in \{f_0, f_1\}} \bar{\mathbb{P}}_{f_j}^n(\psi(D) \neq j)}_{p_{e,1}}\]
Now, by Tsybakov Thm 2.2 (iii), if $KL(\bar{\mathbb{P}}_0^n, \bar{\mathbb{P}}_1^n) \le \beta$, then $p_{e,1} \ge \max \left(\frac{1}{4} \exp(-\beta), \frac{1-\sqrt{\beta/2}}{2}\right)$. 
Now, by \Cref{lem:KL_neff} we know that $KL(\bar{\mathbb{P}}_0^n, \bar{\mathbb{P}}_1^n) \le n_{\eff}KL(\mathbb{P}_0, \mathbb{P}_1)$.
In particular, for Hölder-$\Sigma(\beta,L)$ kernel-bump alternatives of bandwidth $h\to0$,
\[
\KL\!\left(\bar P_{\theta_1}^{(n)}\,\big\|\,\bar P_{\theta_0}^{(n)}\right)
= n_{\eff}\, c\,L^2 h^{2\beta+1}+o(1),
\]
for a constant $c>0$ depending on the kernel and noise law. Thus, choosing $h = n_\eff^{-\frac{1}{2\beta + 1}} \left(\frac{\beta}{K_{\max}^2 L^2 c}\right)^{\frac{1}{2\beta + 1}}$, guarantees $KL(\bar{\mathbb{P}}_0^n, \bar{\mathbb{P}}_1^n) \le \beta$ and thus
\[\inf_{\hat{f}_n} \sup_{f \in \Sigma(\beta, L)}\E_\Xi \E_{\xi, f, n}[(\hat{f}_n(x_0) - f(x_0))^2] \ge n_{\eff}^{-\frac{2\beta}{2\beta + 1}} c_\beta\]
where $c_\beta$ is a constant that depends on $L$, $\beta$ and $\sigma^2$, as we can choose the kernel to be bounded by 1. Taking limits as $n \to \infty$ completes the proof. 
\subsection{Proof of \Cref{lem:KL_neff}}
\label{subsec:kl_scaling_proof}
\begin{proof}[Proof of \Cref{lem:KL_neff}]
    
Recall the correlated noise model. We split the domain of $X$ into $B_X$ buckets and parameterize the shift by $\xi =\{\xi_b\}_{b=1}^{B_X}$, with $\xi_b \sim^{i.i.d.} \mathcal{N}(0, \delta^2 \sigma^2)$. 
\[P_\xi(x, \e) = P_0^X(x)P_0^\e(\e - \xi_{b(x)}),\]
where $b(x)$ is the bucket for $x$. Let $m_b = |\{X \in \{X_1, \dots, X_n\}|X \in b\}|$ be the number of points in bucket $b$ and let $m = \frac{n}{B_X}$ denote the average. By \Cref{ass:0noise_normal}, when data is drawn from the mixture model $\bar{P}^n = \int_{\Xi} P^{\xi, n} d\xi$, $(X_i, Y_i)_{i=1}^n \sim \bar{P}^n$ actually follows a multivariate gaussian distribution, where
\[Y_{1:n} \sim \mathcal{N}\left(f(X_{1:n}), \Sigma\right), \quad \Sigma = \sigma^2I_{B_X} \otimes \Sigma_{\mathrm{blk}}
\quad(\in\R^{n\times n}), \quad \Sigma_{\mathrm{blk}, b}= I_{m_b} + \delta^2\1_{m_b}\1_{m_b}^\top
\quad(\in\R^{m_b \times m_b})\]
Thus, the within-bucket $m_b \times m_b$ matrix has diagonal entries $(1 + \delta^2)\sigma^2$ and off-diagonal entries $\delta^2\sigma^2$. By Sherman-Morrison we have that the inverse matrices are \[\Sigma_{blk,b}^{-1} = \frac{1}{\sigma^2}\left(I_{m_b} - \frac{\delta^2}{1 + m_b \delta^2}\1_{m_b}\1_{m_b}^\top\right)\]
Since the $X$ distribution is fixed, $KL(\bar{P}_0^n||\bar{P}_1^n) = \E_X[KL(\bar{P}_0^n(Y_{1:n}|X_{1:n})||\bar{P}_1^n(Y_{1:n}|X)_{1:n})]$. Under $\bar{P}_j$, $Y_{1:n}|X_{1:n} \sim \mathcal{N}(f_j(X_{1:n}), \Sigma)$, so  

\begin{align*}
    KL(\bar{P}_0^n(Y_{1:n}|X_{1:n})||\bar{P}_1^n(Y_{1:n}|X)_{1:n}) &= \frac{1}{2}(f_1(X_{1:n}) - f_0(X_{1:n}))^T \Sigma^{-1}(f_1(X_{1:n}) - f_0(X_{1:n}))\\
    (f_0 \equiv 0)\quad &= \frac{1}{2\sigma^2} \sum_{b=1}^{B_X} \left(f_1(X_{1:n})^T_b \Sigma^{-1}_{blk,b}f_1(X_{1:n})_b \right) \\
    &= \frac{1}{2 \sigma^2} \sum_{b=1}^{B_X} \left(\sum_{i \in b} f_1(X_i)^2 - \frac{\delta^2}{1 + m_b \delta^2 } \left(\sum_{i \in b} f_1(X_i) \right)^2 \right)
    \intertext{Now set $\hat{f}_b = \frac{1}{m_b} \sum_{i \in b} f_1(X_i)$, and for $X_i=x$, set $d_i = f_1(X_i) - \hat{f}_{b(x)}$. Note that for all $b$, $\sum_{i \in b} d_i = 0$.}
    &= \frac{1}{2\sigma^2} \sum_{b=1}^{B_X} \left(\sum_{i \in b} (\hat{f}_b + d_i)^2 - \frac{\delta^2}{1 + m_b \delta^2} \left(m_b\hat{f}_b\right)^2\right)\\
    &= \frac{1}{2\sigma^2} \sum_{b=1}^{B_X} \left(m_b \hat{f}_b^2 + \sum_{i \in b} d_i^2- \frac{\delta^2}{1 + m_b \delta^2}  m_b^2 \hat{f}_b^2 \right)
    \end{align*}
    Now taking the expectation over the distribution of $X_{1:n}$, 
    \begin{align*}
    KL(\bar{P}_0^n||\bar{P}_1^n) &= \E_X\left[\frac{1}{2\sigma^2} \sum_{b=1}^{B_X} \left(m_b \hat{f}_b^2 + \sum_{i \in b} d_i^2- \frac{\delta^2}{1 + m_b \delta^2}  m_b^2 \hat{f}_b^2 \right)\right]\\
    \text{(Jensen's) } \quad &\le \frac{1}{2\sigma^2} m \left( 1 - \frac{m \delta^2}{1 + m \delta^2}\right) \E_X\left[ \sum_{b=1}^{B_X}\hat{f}_b^2 \right]+ \E_X\left[ \sum_{i=1}^n d_i^2\right]\\
    &= \frac{1}{2\sigma^2} m \left(\frac{1}{1 + m \delta^2}\right) \E_X\left[ \sum_{b=1}^{B_X}\hat{f}_b^2 \right]+ \E_X\left[ \sum_{i=1}^n d_i^2\right]\\
    &= \frac{1}{2\sigma^2} \left(\frac{n}{1 + \frac{n}{B_X}\delta^2}\right) \frac{1}{B_X} \E_X\left[ \sum_{b=1}^{B_X}\hat{f}_b^2 \right]+ \E_X\left[ \sum_{i=1}^n d_i^2\right]
\end{align*}
Now, setting $f_b = \E[f(X)|X \in b]$, we have
\[\E\left[\hat{f}_b^2\right] \le \E\left[\hat{f}_b^2|m_b > 0\right] = f_b^2 + \Var(f(X)|X \in b)\mathbb{E}\left[\frac{1}{m_b}|m_b > 0\right] \le f_b^2 + \Var(f(X)|X \in b)\]
Since $f_1(x) = Lh^{\beta}K \left(\frac{x-x_0}{h}\right)$, it is $L_k h^{\beta - 1}$ Lipschitz, where $L_k = L||\nabla K||_\infty$. Combining this with the fact that for $X$ supported on $[a,b]$, $\Var(X) \le \frac{1}{4}(b-a)^2$, and plugging it in above, 
\[\E\left[\hat{f}_b^2\right] \le \left(\frac{1}{4}\left(\frac{L_k h^{\beta - 1}}{B_X}\right)^2 + (Lh^{\beta}||K||_\infty)^2\right) \1 \{b \cap [x_0 -h/2, x_0 + h/2] \neq \emptyset\}\]
Where the indicator reflects that $K$ vanishes outside of $[-1/2,1/2]$. Since an $h$ fraction of the buckets have nonzero $f$ values 
\[\frac{1}{B_X}\sum_{b = 1}^{B_X} \E\left[\hat{f}_b^2\right] \le h \times h^{2\beta} \left(\frac{L_k^2}{4B_X^2 h^2} + L^2 K_{\max}^2\right)  = h^{2\beta + 1}CL^2 K_{\max}^2 + O\left(\frac{h^{2\beta -1}}{B_X^2}\right)\]
Again, since $f$ is $L_k h^{\beta-1}$-Lipschitz, we know $d_i \le \frac{L_k h^{\beta-1}}{B_X} \1 \{b \cap [x_0 -h/2, x_0 + h/2] \neq \emptyset\}$
\[\sum_{i=1}^n d_i^2 = \sum_{b=1}^{B_X}\sum_{i \in b}d_i^2 \le \sum_{b=1}^{B_X} \left(\1 \{b \cap [x_0 -h/2, x_0 + h/2] \neq \emptyset\}\sum_{i \in b}\left(\frac{L_k h^{\beta-1}}{B_X}\right)^2\right)  \le L_k^2 h^{2\beta-1}\frac{n}{B_X^2}\]
Putting everything together we have
\begin{align*}
    KL(\bar{P}_0||\bar{P}_1) &\le \frac{1}{2\sigma^2}\left(\frac{n}{1 + n \frac{\delta^2}{B_X}}\right) \left(C L^2 K^2_{\max} h^{2\beta + 1} + O(h^{2\beta - 1}/(B_X^2))\right) + L_k^2 h^{2\beta-1}\frac{n}{B_X^2}\\
    &= C_{\sigma, L,K} n_{\eff}h^{2\beta + 1} + O\left(\frac{h^{2\beta - 1}n}{B_X^2}\right)
\end{align*}
We set $h = n_{\eff}^{-1/(2\beta +1)}$ to get that the first term is constant and the second term behaves as

\[n_{\eff}^{-\frac{2\beta - 1}{2\beta +1}}nB_X^{-2} \approx n^{1-\frac{2\beta - 1}{2\beta +1}}B_X^{-2}(1 + nB_X)^{\frac{2\beta - 1}{2\beta +1}} \approx \max(n^{1-\frac{2\beta - 1}{2\beta +1}}B_X^{-2}, nB_X^{-2 +\frac{2\beta - 1}{2\beta +1}}) = \max(n^{\frac{2}{2\beta + 1}}B_X^{-2}, nB_X^{-\frac{2\beta + 3}{2\beta + 1}})\]
Since $2 > \frac{2}{2\beta + 1}$, and $\frac{2\beta + 3}{2\beta +1 } > 1$ for any $\beta >0$, in the regime where $n/{B_X} \to c$, these terms are $o(1)$.

\end{proof}
\end{document}